\documentclass{article}
\usepackage{amssymb}
\usepackage{amsmath}
\usepackage{harvard}
\usepackage{latexsym}
\usepackage{amsmath,amsfonts,amssymb}
\usepackage{ifthen}
\usepackage{harvard}
\usepackage{latexsym}
\usepackage[nohead]{geometry}

\provideboolean{uglystyle}
\setboolean{uglystyle}{true}
\newtheorem{theorem}{Theorem}[section]

\newtheorem{lemma}[theorem]{Lemma}
\newtheorem{proposition}[theorem]{Proposition}

\newtheorem{remark}[theorem]{Remark}
\newcounter{tcnt}[theorem]

\newcounter{pcnt}[theorem]

\newcounter{ccnt}[theorem]

\newcounter{rcnt}[theorem]

\ifthenelse{\boolean{uglystyle}}{
\sloppy

\usepackage{leqno}
}{
}
\input{tcilatex}

\fontsize{13pt}{13pt} \selectfont
\renewcommand{\O}{{\cal O}}

\newcommand{\myenumi}{\renewcommand{\theenumi}{\alph{enumi}}}
\myenumi

\begin{document}

\title{\textsc{Can One Estimate The Unconditional Distribution of
Post-Model-Selection Estimators? }}
\author{Hannes Leeb \\
Department of Statistics, Yale University\\
and \and Benedikt M. P\"{o}tscher \\
Department of Statistics, University of Vienna}
\date{First version: \ April 2005\\
Revised version: \ February 2007}
\maketitle

\begin{abstract}
We consider the problem of estimating the unconditional distribution of a
post-model-selection estimator. The notion of a post-model-selection
estimator here refers to the combined procedure resulting from first
selecting a model (e.g., by a model selection criterion like AIC or by a
hypothesis testing procedure) and then estimating the parameters in the
selected model (e.g., by least-squares or maximum likelihood), all based on
the same data set. We show that it is impossible to estimate the
unconditional distribution with reasonable accuracy even asymptotically. In
particular, we show that no estimator for this distribution can be uniformly
consistent (not even locally). This follows as a corollary to (local)
minimax lower bounds on the performance of estimators for the distribution;
performance is here measured by the probability that the estimation error
exceeds a given threshold. These lower bounds are shown to approach $1/2$ or
even $1$ in large samples, depending on the situation considered. Similar
impossibility results are also obtained for the distribution of linear
functions (e.g., predictors) of the post-model-selection estimator.

{\footnotesize \textsc{AMS Mathematics Subject Classification 2000:} 62F10,
62F12, 62J05, 62J07, 62C05.\vspace{0.2cm} \newline
}

{\footnotesize \ Research of the first author was supported by the Max Kade
Foundation and by the Austrian National Science Foundation (FWF), Grant No.
P13868-MAT. A preliminary draft of the material in this paper was already
written in 1999. }
\end{abstract}

\section{Introduction and Overview\label{s1}}

In many statistical applications a data-based model selection step precedes
the final parameter estimation and inference stage. For example, the
specification of the model (choice of functional form, choice of regressors,
number of lags, etc.) is often based on the data. In contrast, the
traditional theory of statistical inference is concerned with the properties
of estimators and inference procedures under the central assumption of an a
priori given model. That is, it is assumed that the model is known to the
researcher prior to the statistical analysis, except for the value of the
true parameter vector. As a consequence, the actual statistical properties
of estimators or inference procedures following a data-driven model
selection step are not described by the traditional theory which assumes an
a priori given model; in fact, they may differ substantially from the
properties predicted by this theory, cf., e.g., Danilov and Magnus (2004),
Dijkstra and Veldkamp (1988), P\"{o}tscher (1991, Section 3.3), or Rao and
Wu (2001, Section 12). Ignoring the additional uncertainty originating from
the data-driven model selection step and (inappropriately) applying
traditional theory can hence result in very misleading conclusions.

Investigations into the distributional properties of post-model-selection
estimators, i.e., of estimators constructed after a data-driven model
selection step, are relatively few and of recent vintage. Sen (1979)
obtained the unconditional large-sample limit distribution of a
post-model-selection estimator in an i.i.d. maximum likelihood framework,
when selection is between two competing nested models. In P\"{o}tscher
(1991) the asymptotic properties of a class of post-model-selection
estimators (based on a sequence of hypothesis tests) were studied in a
rather general setting covering non-linear models, dependent processes, and
more than two competing models. In that paper, the large-sample limit
distribution of the post-model-selection estimator was derived, both
unconditional as well as conditional on having chosen a correct model, not
necessarily the minimal one. See also P\"{o}tscher and Novak (1998) for
further discussion and a simulation study, and Nickl (2003) for extensions.
The finite-sample distribution of a post-model-selection estimator, both
unconditional and conditional on having chosen a particular (possibly
incorrect) model, was derived in Leeb and P\"{o}tscher (2003) in a normal
linear regression framework; this paper also studied asymptotic
approximations that are in a certain sense superior to the asymptotic
distribution derived in P\"{o}tscher (1991). The distributions of
corresponding linear predictors constructed after model selection were
studied in Leeb (2005, 2006). Related work can also be found in Sen and
Saleh (1987), Kabaila (1995), P\"{o}tscher (1995), Ahmed and Basu (2000),
Kapetanios (2001), Hjort and Claeskens (2003), Duki\'{c} and Pe\~{n}a
(2005), and Leeb and P\"{o}tscher (2005a). The latter paper provides a
simple exposition of the problems of inference post model selection and may
serve as an entry point to the present paper.

It transpires from the papers mentioned above that the finite-sample
distributions (as well as the large-sample limit distributions) of
post-model-selection estimators typically depend on the unknown model
parameters, often in a complicated fashion. For inference purposes, e.g.,
for the construction of confidence sets, estimators of these distributions
would be desirable. Consistent estimators of these distributions can
typically be constructed quite easily, e.g., by suitably replacing unknown
parameters in the large-sample limit distributions by estimators; cf.
Section \ref{s4.1}. However, the merits of such `plug-in' estimators in
small samples are questionable: It is known that the convergence of the
finite-sample distributions to their large-sample limits is typically not
uniform with respect to the underlying parameters (see Appendix \ref{ae}
below and Corollary 5.5 in Leeb and P\"{o}tscher (2003)), and there is no
reason to believe that this non-uniformity will disappear when unknown
parameters in the large-sample limit are replaced by estimators. This
observation is the main motivation for the present paper to investigate in
general the performance of estimators of the distribution of a
post-model-selection estimator, where the estimators of the distribution are
not necessarily `plug-in' estimators based on the limiting distribution. In
particular, we ask whether estimators of the distribution function of
post-model-selection estimators exist that do not suffer from the
non-uniformity phenomenon mentioned above. As we show in this paper the
answer in general is `No'. We also show that these negative results extend
to the problem of estimating the distribution of linear functions (e.g.,
linear predictors) of post-model-selection estimators. Similar negative
results apply also to the estimation of the mean squared error or bias of
post-model-selection estimators; cf. Remark \ref{r6.8}.

To fix ideas consider for the moment the linear regression model 
\begin{equation}
Y=V\chi +W\psi +u  \label{1.1}
\end{equation}%
where $V$ and $W$, respectively, represent $n\times k$ and $n\times l$
non-stochastic regressor matrices ($k\geq 1,l\geq 1$), and the $n\times 1$
disturbance vector $u$ is normally distributed with mean zero and
variance-covariance matrix $\sigma ^{2}I_{n}$. We also assume for the moment
that $(V:W)^{\prime }(V:W)/n$ converges to a non-singular matrix as the
sample size $n$ goes to infinity and that $\lim_{n\rightarrow \infty
}V^{\prime }W/n\neq 0$ (for a discussion of the case where this limit is
zero see Example 1 in Section \ref{s5}\texttt{)}. Now suppose that the
vector $\chi $ represents the parameters of interest, while the parameter
vector $\psi $ and the associated regressors in $W$ have been entered into
the model only to avoid possible misspecification. Suppose further that the
necessity to include $\psi $ or some of its components is then checked on
the basis of the data, i.e., a model selection procedure is used to
determine which components of $\psi $ are to be retained in the model, the
inclusion of $\chi $ not being disputed. The selected model is then used to
obtain the final (post-model-selection) estimator $\tilde{\chi}$ for $\chi $%
. We are now interested in the unconditional finite-sample distribution of $%
\tilde{\chi}$ (appropriately scaled and centered). Denote this $k$%
-dimensional cumulative distribution function (cdf) by $G_{n,\theta ,\sigma
}(t)$. As indicated in the notation, this distribution function depends on
the true parameters $\theta =(\chi ^{\prime },\psi ^{\prime })^{\prime }$
and $\sigma $. For the sake of definiteness of discussion assume for the
moment that the model selection procedure used here is the particular
`general-to-specific' procedure described at the beginning of Section \ref%
{s2}; we comment on other model selection procedures, including Akaike's AIC
and thresholding procedures, below.

As mentioned above, it is not difficult to construct a consistent estimator
of $G_{n,\theta ,\sigma }(t)$ for any $t$, i.e., an estimator $\hat{G}%
_{n}(t) $ satisfying 
\begin{equation}
P_{n,\theta ,\sigma }\left( \left\vert \hat{G}_{n}(t)-G_{n,\theta ,\sigma
}(t)\right\vert \,>\,\delta \right) \overset{n\rightarrow \infty }{%
\longrightarrow }0  \label{1.2}
\end{equation}%
for each $\delta >0$ and each $\theta $, $\sigma $; see Section \ref{s4.1}.
However, it follows from the results in Section \ref{s5} that \emph{any }%
estimator satisfying (\ref{1.2}), i.e., \emph{any consistent }estimator of $%
G_{n,\theta ,\sigma }(t)$, necessarily also satisfies 
\begin{equation}
\sup_{\,||\theta ||<R}\;P_{n,\theta ,\sigma }\left( \left\vert \hat{G}%
_{n}(t)-G_{n,\theta ,\sigma }(t)\right\vert \,>\,\delta \right) \overset{%
n\rightarrow \infty }{\longrightarrow }1  \label{1.3}
\end{equation}%
for suitable positive constants $R$ and $\delta $ that do not depend on the
estimator. That is, while the probability in (\ref{1.2}) converges to zero
for every given $\theta $ by consistency, relation (\ref{1.3}) shows that it
does not do so uniformly in $\theta $. It follows that $\hat{G}_{n}(t)$ can
never be uniformly consistent (not even when restricting consideration to
uniform consistency over all compact subsets of the parameter space). Hence,
a large sample size does not guarantee a small estimation error with high
probability when estimating the distribution function of a
post-model-selection estimator. In this sense, reliably assessing the
precision of post-model-selection estimators is an intrinsically hard
problem. Apart from (\ref{1.3}), we also provide minimax lower bounds for
arbitrary (not necessarily consistent) estimators of the conditional
distribution function $G_{n,\theta ,\sigma }(t)$. For example, we provide
results that imply that 
\begin{equation}
\liminf_{n\rightarrow \infty }\inf_{\hat{G}_{n}(t)}\sup_{\,||\theta
||<R}\;P_{n,\theta ,\sigma }\left( \left\vert \hat{G}_{n}(t)-G_{n,\theta
,\sigma }(t)\right\vert \,>\,\delta \right) >0  \label{1.4}
\end{equation}%
holds for suitable positive constants $R$ and $\delta $, where the infimum
extends over \emph{all} estimators of $G_{n,\theta ,\sigma }(t)$. The
results in Section \ref{s5} in fact show that the balls $||\theta ||<R$ in (%
\ref{1.3}) and (\ref{1.4}) can be replaced by suitable balls (not
necessarily centered at the origin) shrinking at the rate $n^{-1/2}$. This
shows that the non-uniformity phenomenon described in (\ref{1.3})-(\ref{1.4}%
) is a local, rather than a global, phenomenon. In Section \ref{s5} we
further show that the non-uniformity phenomenon expressed in (\ref{1.3}) and
(\ref{1.4}) typically also arises when the parameter of interest is not $%
\chi $, but some other linear transformation of $\theta =(\chi ^{\prime
},\psi ^{\prime })^{\prime }$. As discussed in Remark \ref{r6.1}, the
results also hold for randomized estimators of the unconditional
distribution function $G_{n,\theta ,\sigma }(t)$. Hence no resampling
procedure whatsoever can alleviate the problem. This explains the anecdotal
evidence in the literature that resampling methods are often unsuccessful in
approximating distributional properties of post-model-selection estimators
(e.g., Dijkstra and Veldkamp (1988), or Freedman, Navidi, and Peters
(1988)). See also the discussion on resampling in Section \ref{s7}.

The results outlined above are presented in Section \ref{s4} for the
particular `general-to-specific' model selection procedure described at the
beginning of Section \ref{s2}. Analogous results for a large class of model
selection procedures, including Akaike's AIC and thresholding procedures,
are then given in Section \ref{s99}, based on the results in Section \ref{s4}%
. In fact, the non-uniformity phenomenon expressed in (\ref{1.3})-(\ref{1.4}%
) is not specific to the model selection procedures discussed in Sections %
\ref{s2} and \ref{s99} of the present paper, but will occur for most (if not
all) model selection procedures, including consistent ones; cf. Sections \ref%
{after}\ and \ref{s7} for more discussion. Section \ref{after} also shows
that the results are -- as is to be expected -- by no means limited to the
linear regression model.

We focus on the unconditional distributions of post-model-selection
estimators in the present paper. One can, however, also envisage a situation
where one is more interested in the conditional distribution given the
outcome of the model selection procedure. In line with the literature on
conditional inference (see, e.g., Robinson (1979) or Lehmann and Casella
(1998, p. 421)), one may argue that, given the outcome of the model
selection step, the relevant object of interest is the conditional rather
than the unconditional distribution of the post-model-selection estimator.
In this case similar results can be obtained and are reported in Leeb and P%
\"{o}tscher (2006b). We note that on a technical level the results in Leeb
and P\"{o}tscher (2006b) and in the present paper require separate treatment.

The plan of the paper is as follows: Post-model-selection estimators based
on a `general-to-specific' model selection procedure are the subject of
Section \ref{s2}. After introducing the basic framework and some notation,
like the family of models $M_{p}$ from which the `general-to-specific' model
selection procedure $\hat{p}$ selects, as well as the post-model-selection
estimator $\tilde{\theta}$, the unconditional cdf $G_{n,\theta ,\sigma }(t)$
of (a linear function of) the post-model-selection estimator $\tilde{\theta}$
is discussed in Section \ref{s3}. Consistent estimators of $G_{n,\theta
,\sigma }(t)$ are given in Section \ref{s4.1}. The main results of the paper
are contained in Section \ref{s5} and Section \ref{s99}: In Section \ref{s5}
we provide a detailed analysis of the non-uniformity phenomenon encountered
in (\ref{1.3})-(\ref{1.4}). In Section \ref{s99} the `impossibility' result
from Section \ref{s5} is extended to a large class of model selection
procedures including Akaike's AIC and to selection from a non-nested
collection of models. Some remarks are collected in Section \ref{s6}, while
Section \ref{after} discusses extensions and the scope of the results of the
paper. Conclusions are drawn in Section \ref{s7}. All proofs as well as some
auxiliary results are collected into appendices. Finally a word on notation:
The Euclidean norm is denoted by $\left\Vert \cdot \right\Vert $, and $%
\lambda _{\max }(E)$ denotes the largest eigenvalue of a symmetric matrix $E$%
. A prime denotes transposition of a matrix. For vectors $x$ and $y$ the
relation $x\leq y$ ($x<y$, respectively) denotes $x_{i}\leq y_{i}$ ($%
x_{i}<y_{i}$, respectively) for all $i$. As usual, $\Phi $ denotes the
standard normal distribution function.

\section{Results for Post-Model-Selection Estimators Based on a
`General-to-Specific' Model Selection Procedure \label{s2}}

Consider the linear regression model 
\begin{equation}
Y=X\theta +u,  \label{2.1}
\end{equation}%
where $X$ is a non-stochastic $n\times P$ matrix with $rank(X)=P$ and $u\sim
N(0,\sigma ^{2}I_{n})$, $\sigma ^{2}>0$. Here $n$ denotes the sample size
and we assume $n>P\geq 1$. In addition, we assume that $Q=\lim_{n\rightarrow
\infty }X^{\prime }X/n$ exists and is non-singular. In this section we shall
-- similar as in P\"{o}tscher (1991) -- consider model selection from the
collection of nested models $M_{\O }\subseteq M_{\O +1}\subseteq \dots
\subseteq M_{P}$, where $\O $ is specified by the user, and where for $0\leq
p\leq P$ the model $M_{p}$\ is given by 
\begin{equation*}
M_{p}=\left\{ (\theta _{1},\dots ,\theta _{P})^{\prime }\in \mathbf{R}%
^{P}:\,\theta _{p+1}=\dots =\theta _{P}=0\right\} .
\end{equation*}%
[In Section \ref{s99} below also general non-nested families of models will
be considered.] Clearly, the model $M_{p}$ corresponds to the situation
where only the first $p$ regressors in (\ref{2.1}) are included. For the
most parsimonious model under consideration, i.e., for $M_{\O }$, we assume
that $\O $ satisfies $0\leq \O <P$; if $\O >0$, this model contains as free
parameters only those components of the parameter vector $\theta $ that are
not subject to model selection. [In the notation used in connection with (%
\ref{1.1}) we then have $\chi =(\theta _{1},\dots ,\theta _{\O })^{\prime }$
and $\psi =$ $(\theta _{\O +1},\dots ,\theta _{P})^{\prime }$.] Furthermore,
note that $M_{0}=\{(0,\dots ,0)^{\prime }\}$ and that $M_{P}=\mathbf{R}^{P}$%
. We call $M_{p}$ the regression model of order $p$.

The following notation will prove useful. For matrices $B$ and $C$ of the
same row-dimension, the column-wise concatenation of $B$ and $C$ is denoted
by $(B:C)$. If $D$ is an $m\times P$ matrix, let $D[p]$ denote the $m\times
p $ matrix consisting of the first $p$ columns of $D$. Similarly, let $%
D[\lnot p]$ denote the $m\times (P-p)$ matrix consisting of the last $P-p$
columns of $D$. If $x$ is a $P\times 1$ vector, we write in abuse of
notation $x[p]$ and $x[\lnot p]$ for $(x^{\prime }[p])^{\prime }$ and $%
(x^{\prime }[\lnot p])^{\prime }$, respectively. [We shall use the above
notation also in the `boundary' cases $p=0$ and $p=P$. It will always be
clear from the context how expressions containing symbols like $D[0]$, $%
D[\lnot P]$, $x[0]$, or $x[\lnot P]$ are to be interpreted.] As usual, the $%
i $-th component of a vector $x$ is denoted by $x_{i}$, and the entry in the 
$i $-th row and $j$-th column of a matrix $B$ is denoted by $B_{i,j}$.

The restricted least-squares estimator of $\theta $ under the restriction $%
\theta \lbrack \lnot p]=0$, i.e., under $\theta _{p+1}=\dots =\theta _{P}=0$%
, will be denoted by $\tilde{\theta}(p)$, $0\leq p\leq P$ (in case $p=P$ the
restriction being void). Note that $\tilde{\theta}(p)$ is given by the $%
P\times 1$ vector 
\begin{equation}
\tilde{\theta}(p)=\left( 
\begin{array}{c}
\left( X[p]^{\prime }X[p]\right) ^{-1}X[p]^{\prime }Y \\ 
(0,\dots ,0)^{\prime }%
\end{array}%
\right) ,  \notag
\end{equation}%
where the expressions $\tilde{\theta}(0)$ and $\tilde{\theta}(P)$,
respectively, are to be interpreted as the zero-vector in $\mathbf{R}^{P}$
and as the unrestricted least-squares estimator of $\theta $. Given a
parameter vector $\theta $ in $\mathbf{R}^{P}$, the order of $\theta $
(relative to the nested sequence of models $M_{p}$) is defined as 
\begin{equation}
p_{0}(\theta )=\min \left\{ p:\;0\leq p\leq P,\;\theta \in M_{p}\right\} . 
\notag  \label{2'}
\end{equation}%
Hence, if $\theta $ is the true parameter vector, a model $M_{p}$ is a
correct model if and only if $p\geq p_{0}(\theta )$. We stress that $%
p_{0}(\theta )$ is a property of a \emph{single parameter}, and hence needs
to be distinguished from the notion of the order of the model $M_{p}$
introduced earlier, which is a property of the \emph{set of parameters} $%
M_{p}$.

A model selection procedure is now nothing else than a data-driven
(measurable) rule $\hat{p}$ that selects a value from $\{\O ,\dots ,P\}$ and
thus selects a model from the list of candidate models $M_{\O },\dots ,M_{P}$%
. In this section we shall consider as an important leading case a
`general-to-specific' model selection procedure based on a sequence of
hypothesis tests. [Results for a larger class of model selection procedures,
including Akaike's AIC, are provided in Section \ref{s99}.] This procedure
is given as follows: The sequence of hypotheses $H_{0}^{p}:\,p_{0}(\theta
)<p $ is tested against the alternatives $H_{1}^{p}:\,p_{0}(\theta )=p$ in
decreasing order starting at $p=P$. If, for some $p>\O $, $H_{0}^{p}$ is the
first hypothesis in the process that is rejected, we set $\hat{p}=p$. If no
rejection occurs until even $H_{0}^{\O +1}$ is not rejected, we set $\hat{p}=%
\O $. Each hypothesis in this sequence is tested by a kind of $t$-test where
the error variance is always estimated from the overall model (but see the
discussion following Theorem \ref{t99.2} in Section \ref{s99} below for
other choices of estimators of the error variance). More formally, we have 
\begin{equation}
\hat{p}=\max \left\{ p:\,|T_{p}|\geq c_{p},\;0\leq p\leq P\right\} ,
\label{2.2}
\end{equation}%
with $c_{\O }=0$ in order to ensure a well-defined $\hat{p}$ in the range $\{%
\O ,\O +1,\dots ,P\}$. For $\O <p\leq P$, the critical values $c_{p}$
satisfy $0<c_{p}<\infty $ and are independent of sample size (but see also
Remark \ref{r6.4}). The test-statistics are given by 
\begin{equation}
T_{p}=\frac{\sqrt{n}\tilde{\theta}_{p}(p)}{\hat{\sigma}\xi _{n,p}}\qquad
(0<p\leq P)  \notag  \label{4}
\end{equation}%
with the convention that $T_{0}=0$. Furthermore, 
\begin{equation}
\xi _{n,p}=\left( \left[ \left( \frac{X[p]^{\prime }X[p]}{n}\right) ^{-1}%
\right] _{p,p}\right) ^{\frac{1}{2}}\qquad (0<p\leq P)  \notag  \label{4'}
\end{equation}%
denotes the nonnegative square root of the $p$-th diagonal element of the
matrix indicated, and $\hat{\sigma}^{2}$ is given by 
\begin{equation}
\hat{\sigma}^{2}=(n-P)^{-1}(Y-X\tilde{\theta}(P))^{\prime }(Y-X\tilde{\theta}%
(P)).  \notag  \label{5}
\end{equation}%
Note that under the hypothesis $H_{0}^{p}$ the statistic $T_{p}$ is $t$%
-distributed with $n-P$ degrees of freedom for $0<p\leq P$. It is also easy
to see that the so-defined model selection procedure $\hat{p}$ is
conservative: The probability of selecting an incorrect model, i.e., the
probability of the event $\{\hat{p}<p_{0}(\theta )\}$, converges to zero as
the sample size increases. In contrast, the probability of the event $\{\hat{%
p}=p\}$, for $p$ satisfying $\max \{p_{0}(\theta ),\O \}\leq p\leq P$,
converges to a positive limit; cf., for example, Proposition~5.4 and
equation (5.6) in Leeb (2006).

The post-model-selection estimator $\tilde{\theta}$ can now be defined as
follows: On the event $\hat{p}=p$, $\tilde{\theta}$ is given by the
restricted least-squares estimator $\tilde{\theta}(p)$, i.e., 
\begin{equation}
\tilde{\theta}=\sum_{p=\O }^{P}\tilde{\theta}(p)\,\mathbf{1(}\hat{p}=p),
\label{2.3}
\end{equation}%
where $\mathbf{1(\cdot )}$ denotes the indicator function of the event shown
in the argument.

\subsection{The Distribution of the Post-Model-Selection Estimator\label{s3}}

We now introduce the distribution function of a linear transformation of $%
\tilde{\theta}$ and summarize some of its properties that will be needed in
the subsequent development. To this end, let $A$ be a non-stochastic $%
k\times P$ matrix of rank $k$, $1\leq k\leq P$, and consider the cdf 
\begin{equation}
G_{n,\theta ,\sigma }(t)=P_{n,\theta ,\sigma }\left( \sqrt{n}A(\tilde{\theta}%
-\theta )\leq t\right) \qquad (t\in \mathbf{R}^{k}).  \label{3.1}
\end{equation}%
Here $P_{n,\theta ,\sigma }(\cdot )$ denotes the probability measure
corresponding to a sample of size $n$ from (\ref{2.1}).

Depending on the choice of the matrix $A$, several important scenarios are
covered by (\ref{3.1}): The cdf of $\sqrt{n}(\tilde{\theta}-\theta )$ is
obtained by setting $A$ equal to the $P\times P$\ identity matrix $I_{P}$.
In case $\O >0$, the cdf of those components of $\sqrt{n}(\tilde{\theta}%
-\theta )$ which correspond to the parameter of interest $\chi $ in (\ref%
{1.1}) can be studied by setting $A$ to the $\O \times P$ matrix $(I_{\O %
}:0) $ as we then have $A\theta =(\theta _{1},\dots ,\theta _{\O })^{\prime
}=\chi $. Finally, if $A\neq 0$ is an $1\times P$ vector, we obtain the
distribution of a linear predictor based on the post-model-selection
estimator. See the examples at the end of Section \ref{s5} for more
discussion.

The cdf $G_{n,\theta ,\sigma }$ and its properties have been analyzed in
detail in Leeb and P\"{o}tscher (2003) and Leeb (2006). To be able to access
these results we need some further notation. Note that on the event $\hat{p}%
=p$ the expression $A(\tilde{\theta}-\theta )$ equals $A(\tilde{\theta}%
(p)-\theta )$ in view of (\ref{2.3}). The expected value of the restricted
least-squares estimator $\tilde{\theta}(p)$ will be denoted by $\eta _{n}(p)$
and is given by the $P\times 1$ vector 
\begin{equation}
\eta _{n}(p)=\left( 
\begin{array}{c}
\theta \lbrack p]\,+\,(X[p]^{\prime }X[p])^{-1}X[p]^{\prime }X[\lnot
p]\theta \lbrack \lnot p] \\ 
(0,\dots ,0)^{\prime }%
\end{array}%
\right)  \label{3.2}
\end{equation}%
with the conventions that $\eta _{n}(0)=(0,\dots ,0)^{\prime }\in \mathbf{R}%
^{P}$ and that $\eta _{n}(P)=\theta $. Furthermore, let $\Phi _{n,p}$ denote
the cdf of $\sqrt{n}A(\tilde{\theta}(p)-\eta _{n}(p))$, i.e., the cdf of $%
\sqrt{n}A$ times the restricted least-squares estimator based on model $%
M_{p} $ centered at its mean. Hence, $\Phi _{n,p}$\ is the cdf of a $k$%
-variate Gaussian random vector with mean zero and variance-covariance
matrix $\sigma ^{2}A[p](X[p]^{\prime }X[p]/n)^{-1}A[p]^{\prime }$ in case $%
p>0$, and it is the cdf of point-mass at zero in $\mathbf{R}^{k}$ in case $%
p=0$. If $p>0$ and if the matrix $A[p]$ has full row rank $k$, then $\Phi
_{n,p}$ has a density with respect to Lebesgue measure, and we shall denote
this density by $\phi _{n,p}$. We note that $\eta _{n}(p)$ depends on $%
\theta $ and that $\Phi _{n,p}$ depends on $\sigma $ (in case $p>0$),
although these dependencies are not shown explicitly in the notation.

For $p>0$ we introduce 
\begin{equation}
b_{n,p}=C_{n}^{(p)^{\prime }}(A[p](X[p]^{\prime }X[p]/n)^{-1}A[p]^{\prime
})^{-},  \label{3.3}
\end{equation}%
and 
\begin{equation}
\zeta _{n,p}^{2}=\xi _{n,p}^{2}-C_{n}^{(p)^{\prime }}(A[p](X[p]^{\prime
}X[p]/n)^{-1}A[p]^{\prime })^{-}C_{n}^{(p)},  \label{3.4}
\end{equation}%
with $\zeta _{n,p}\geq 0$. Here $C_{n}^{(p)}=$ $A[p](X[p]^{\prime
}X[p]/n)^{-1}e_{p}$, where $e_{p}$ denotes the $p$-th standard basis vector
in $\mathbf{R}^{p}$, and $B^{-}$ denotes a generalized inverse of a matrix $%
B $. [Observe that $\zeta _{n,p}^{2}$ is invariant under the choice of the
generalized inverse. The same is not necessarily true for $b_{n,p}$, but is
true for $b_{n,p}z$ for all $z$ in the column-space of $A[p]$. Also note
that (\ref{3.5}) below depends on $b_{n,p}$ only through $b_{n,p}z$ with $z$
in the column-space of $A[p]$.] We observe that the vector of covariances
between $A\tilde{\theta}(p)$ and $\tilde{\theta}_{p}(p)$ is precisely given
by $\sigma ^{2}n^{-1}C_{n}^{(p)}$ (and hence does \emph{not} depend on $%
\theta $). Furthermore, observe that $A\tilde{\theta}(p)$ and $\tilde{\theta}%
_{p}(p)$ are uncorrelated if and only if $\zeta _{n,p}^{2}=\xi _{n,p}^{2}$
if and only if $b_{n,p}z=0$ for all $z$ in the column-space of $A[p]$; cf.
Lemma A.2 in Leeb (2005).

Finally, for a univariate Gaussian random variable $\mathfrak{N}$ with zero
mean and variance $s^{2}$, $s\geq 0$, we write $\Delta _{s}(a,b)$ for $%
\mathbb{P(}|\mathfrak{N}-a|<b)$, $a\in \mathbf{R}\cup \{-\infty ,\infty \}$, 
$b\in \mathbf{R}$. Note that $\Delta _{s}(\cdot ,\cdot )$ is symmetric
around zero in its first argument, and that $\Delta _{s}(-\infty ,b)=\Delta
_{s}(\infty ,b)=0$ holds. In case $s=0$, $\mathfrak{N}$ is to be interpreted
as being equal to zero, hence $a\mapsto \Delta _{0}(a,b)$ reduces to the
indicator function of the interval $\left( -b,b\right) $.

We are now in a position to present the explicit formula for $G_{n,\theta
,\sigma }(t)$ derived in Leeb (2006): 
\begin{align}
G_{n,\theta ,\sigma }(t)& =\Phi _{n,\O }(t-\sqrt{n}A(\eta _{n}(\O )-\theta
))\int_{0}^{\infty }\prod_{q=\O +1}^{P}\Delta _{\sigma \xi _{n,q}}(\sqrt{n}%
\eta _{n,q}(q),sc_{q}\sigma \xi _{n,q})h(s)ds  \notag \\
& +\sum_{p=\O +1}^{P}\int_{z\leq t-\sqrt{n}A(\eta _{n}(p)-\theta )}\,\Big[%
\int_{0}^{\infty }(1-\Delta _{\sigma \zeta _{n,p}}(\sqrt{n}\eta
_{n,p}(p)+b_{n,p}z,sc_{p}\sigma \xi _{n,p}))  \label{3.5} \\
& \qquad \prod_{q=p+1}^{P}\Delta _{\sigma \xi _{n,q}}(\sqrt{n}\eta
_{n,q}(q),sc_{q}\sigma \xi _{n,q})h(s)ds\Big]\Phi _{n,p}(dz).  \notag
\end{align}%
In the above display, $\Phi _{n,p}(dz)$ denotes integration with respect to
the measure induced by the normal cdf $\Phi _{n,p}$ on $\mathbf{R}^{k}$ and $%
h$ denotes the density of $\hat{\sigma}/\sigma $, i.e., $h$ is the density
of $(n-P)^{-1/2}$ times the square-root of a chi-square distributed random
variable with $n-P$ degrees of freedom. The finite-sample distribution of
the post-model-selection estimator given in (\ref{3.5}) is in general not
normal, e.g., it can be bimodal; see Figure 2 in Leeb and P\"{o}tscher
(2005a) or Figure 1 in Leeb (2006). [An exception where (\ref{3.5}) is
normal is the somewhat trivial case where $C_{n}^{(p)}=0$, i.e., where $A%
\tilde{\theta}(p)$ and $\tilde{\theta}_{p}(p)$ are uncorrelated, for $p=\O %
+1,\ldots ,P$; see Leeb (2006, Section 3.3) for more discussion.] We note
for later use that $G_{n,\theta ,\sigma }(t)=\sum_{p=\O }^{P}G_{n,\theta
,\sigma }(t|p)\pi _{n,\theta ,\sigma }(p)$ where $G_{n,\theta ,\sigma }(t|p)$
represents the cdf of $\sqrt{n}A(\tilde{\theta}-\theta )$ conditional on the
event $\{\hat{p}=p\}$ and where $\pi _{n,\theta ,\sigma }(p)$ is the
probability of this event under $P_{n,\theta ,\sigma }$. Note that $\pi
_{n,\theta ,\sigma }(p)$ is always positive for $\O \leq p\leq P$; cf. Leeb
(2006), Section 3.2.

To describe the large-sample limit of $G_{n,\theta ,\sigma }$, some further
notation is necessary. For $p$ satisfying $0<p\leq P$, partition the matrix $%
Q=\lim_{n\rightarrow \infty }X^{\prime }X/n$ as 
\begin{equation*}
Q=\left( 
\begin{array}{cc}
Q[p:p] & Q[p:\lnot p] \\ 
Q[\lnot p:p] & Q[\lnot p:\lnot p]%
\end{array}%
\right) ,
\end{equation*}%
where $Q[p:p]$ is a $p\times p$ matrix. Let $\Phi _{\infty ,p}$ be the cdf
of a $k$-variate Gaussian random vector with mean zero and
variance-covariance matrix $\sigma ^{2}A[p]Q[p:p]^{-1}A[p]^{\prime }$, $%
0<p\leq P$, and let $\Phi _{\infty ,0}$ denote the cdf of point-mass at zero
in $\mathbf{R}^{k}$. Note that $\Phi _{\infty ,p}$ has a Lebesgue density if 
$p>0$ and the matrix $A[p]$ has full row rank $k$; in this case, we denote
the Lebesgue density of $\Phi _{\infty ,p}$ by $\phi _{\infty ,p}$. Finally,
for $p=1,\dots ,P$, define

\begin{equation*}
\xi _{\infty ,p}^{2}=(Q[p:p]^{-1})_{p,p},
\end{equation*}%
\begin{equation}
\zeta _{\infty ,p}^{2}=\xi _{\infty ,p}^{2}-C_{\infty }^{(p)\prime
}(A[p]Q[p:p]^{-1}A[p]^{\prime })^{-}C_{\infty }^{(p)},  \label{zetainfinity}
\end{equation}%
\begin{equation*}
b_{\infty ,p}=C_{\infty }^{(p)\prime }(A[p]Q[p:p]^{-1}A[p]^{\prime })^{-},
\end{equation*}%
where $C_{\infty }^{(p)}=A[p]Q[p:p]^{-1}e_{p}$, with $e_{p}$ denoting the $p$%
-th standard basis vector in $\mathbf{R}^{p}$; furthermore, take $\zeta
_{\infty ,p}$ and $\xi _{\infty ,p}$ as the nonnegative square roots of $%
\zeta _{\infty ,p}^{2}$ and $\xi _{\infty ,p}^{2}$, respectively. As the
notation suggests, $\Phi _{\infty ,p}$ is the large-sample limit of $\Phi
_{n,p}$, and$\ C_{\infty }^{(p)}$, $\xi _{\infty ,p}^{2}$, and $\zeta
_{\infty ,p}^{2}$ are the limits of $C_{n}^{(p)}$, $\xi _{n,p}^{2}$, and $%
\zeta _{n,p}^{2}$, respectively; moreover, $b_{n,p}z$ converges to $%
b_{\infty ,p}z$ for each $z$ in the column-space of $A[p]$. See Lemma A.2 in
Leeb (2005).

The next result describes the large-sample limit of the cdf under local
alternatives to $\theta $ and is taken from Leeb (2006, Corollary 5.6).
Recall that the total variation distance between two cdfs $G$ and $G^{\ast }$
on $\mathbf{R}^{k}$ is defined as $||G-G^{\ast }||_{TV}$ $%
=\sup_{E}|G(E)-G^{\ast }(E)|$, where the supremum is taken over all Borel
sets $E$. Clearly, the relation $|G(t)-G^{\ast }(t)|\leq ||G-G^{\ast
}||_{TV} $ holds for all $t\in \mathbf{R}^{k}$. Thus, if $G$ and $G^{\ast }$
are close with respect to the total variation distance, then $G(t)$ is close
to $G^{\ast }(t)$, uniformly in $t$.

\begin{proposition}
\label{p3.1} Suppose $\theta \in \mathbf{R}^{P}$ and $\gamma \in \mathbf{R}%
^{P}$ and let $\sigma ^{(n)}$ be a sequence of positive real numbers which
converges to a (finite) limit $\sigma >0$ as $n\rightarrow \infty $. Then
the cdf $G_{n,\theta +\gamma /\sqrt{n},\sigma ^{(n)}}$ converges to a limit $%
G_{\infty ,\theta ,\sigma ,\gamma }$ in total variation, i.e., 
\begin{equation}
\left\vert \left\vert G_{n,\theta +\gamma /\sqrt{n},\sigma ^{(n)}}-G_{\infty
,\theta ,\sigma ,\gamma }\right\vert \right\vert _{TV}\overset{n\rightarrow
\infty }{\longrightarrow }0.  \label{3.6}
\end{equation}%
The large-sample limit cdf $G_{\infty ,\theta ,\sigma ,\gamma }(t)$\ is
given by 
\begin{align}
& \Phi _{\infty ,p_{\ast }}(t-\beta ^{(p_{\ast })})\prod_{q=p_{\ast
}+1}^{P}\Delta _{\sigma \xi _{\infty ,q}}(\nu _{q},c_{q}\sigma \xi _{\infty
,q})  \notag \\
& +\sum_{p=p_{\ast }+1}^{P}\int_{z\leq t-\beta ^{(p)}}(1-\Delta _{\sigma
\zeta _{\infty ,p}}(\nu _{p}+b_{\infty ,p}z,c_{p}\sigma \xi _{\infty
,p}))\Phi _{\infty ,p}(dz)\prod_{q=p+1}^{P}\Delta _{\sigma \xi _{\infty
,q}}(\nu _{q},c_{q}\sigma \xi _{\infty ,q})  \label{3.7}
\end{align}%
where $p_{\ast }=\max \{p_{0}(\theta ),\O \}$. Here for $0\leq p\leq P$%
\begin{equation*}
\beta ^{(p)}=A\left( 
\begin{array}{c}
Q[p:p]^{-1}Q[p:\lnot p]\gamma \lbrack \lnot p] \\ 
-\gamma \lbrack \lnot p]%
\end{array}%
\right) ,
\end{equation*}%
with the convention that $\beta ^{(p)}=-A\gamma $ if $p=0$ and that $\beta
^{(p)}=(0,\ldots ,0)^{\prime }$ if $p=P$. Furthermore, we have set $\nu
_{p}=\gamma _{p}+(Q[p:p]^{-1}Q[p:\lnot p]\gamma \lbrack \lnot p])_{p}$ for $%
p>0$. [Note that $\beta ^{(p)}=\lim_{n\rightarrow \infty }\sqrt{n}A(\eta
_{n}(p)-\theta -\gamma /\sqrt{n})$ for $p\geq p_{0}(\theta )$, and that $\nu
_{p}=\lim_{n\rightarrow \infty }\sqrt{n}\eta _{n,p}(p)$ for $p>p_{0}(\theta
) $. Here $\eta _{n}(p)$ is defined as in (\ref{3.2}), but with $\theta
+\gamma /\sqrt{n}$ replacing $\theta $.]

If $p_{\ast }>0$ and if the matrix $A[p_{\ast }]$ has full row rank $k$,
then the Lebesgue density $\phi _{\infty ,p}$ of $\Phi _{\infty ,p}$ exists
for all $p\geq p_{\ast }$ and hence the density of (\ref{3.7}) exists and is
given by 
\begin{eqnarray*}
&&\phi _{\infty ,p_{\ast }}(t-\beta ^{(p_{\ast })})\prod_{q=p_{\ast
}+1}^{P}\Delta _{\sigma \xi _{\infty ,q}}(\nu _{q},c_{q}\sigma \xi _{\infty
,q}) \\
&&+\sum_{p=p_{\ast }+1}^{P}(1-\Delta _{\sigma \zeta _{\infty ,p}}(\nu
_{p}+b_{\infty ,p}(t-\beta ^{(p)}),c_{p}\sigma \xi _{\infty ,p}))\phi
_{\infty ,p}(t-\beta ^{(p)})\prod_{q=p+1}^{P}\Delta _{\sigma \xi _{\infty
,q}}(\nu _{q},c_{q}\sigma \xi _{\infty ,q}).
\end{eqnarray*}
\end{proposition}

Like the finite-sample distribution, the limiting distribution of the
post-model-selection estimator given in (\ref{3.7}) is in general not
normal. An exception is the case where $C_{\infty }^{(p)}=0$ for $p>p_{\ast
} $ in which case $G_{\infty ,\theta ,\sigma ,\gamma }$ reduces to $\Phi
_{\infty ,P}$;\ see Remark \ref{rb.3} in Appendix \ref{aa}. If $\gamma =0$,
we write $G_{\infty ,\theta ,\sigma }(t)$ as shorthand for $G_{\infty
,\theta ,\sigma ,0}(t)$ in the following.

\subsection{Estimators of the Finite-Sample Distribution\label{s4}}

For the purpose of inference after model selection the finite-sample
distribution of the post-model-selection-estimator is an object of
particular interest. As we have seen, it depends on unknown parameters in a
complicated manner, and hence one will have to be satisfied with estimators
of this cdf. As we shall see, it is not difficult to construct consistent
estimators of $G_{n,\theta ,\sigma }(t)$. However, despite this consistency
result, we shall find in Section \ref{s5}\ that \emph{any} estimator of $%
G_{n,\theta ,\sigma }(t)$ typically performs unsatisfactory, in that the
estimation error can not become small uniformly over (subsets of) the
parameter space even as sample size goes to infinity. In particular, no
uniformly consistent estimators exist, not even locally.

\subsubsection{Consistent Estimators\label{s4.1}}

We construct a consistent estimator of $G_{n,\theta ,\sigma }(t)$ by
commencing from the asymptotic distribution. Specializing to the case $%
\gamma =0$ and $\sigma ^{(n)}=\sigma $\ in Proposition~\ref{p3.1}, the
large-sample limit of $G_{n,\theta ,\sigma }(t)$ is given by 
\begin{align}
G_{\infty ,\theta ,\sigma }(t)& =\Phi _{\infty ,p_{\ast
}}(t)\prod_{q=p_{\ast }+1}^{P}\Delta _{\sigma \xi _{\infty
,q}}(0,c_{q}\sigma \xi _{\infty ,q})  \notag \\
& +\sum_{p=p_{\ast }+1}^{P}\int_{z\leq t}(1-\Delta _{\sigma \zeta _{\infty
,p}}(b_{\infty ,p}z,c_{p}\sigma \xi _{\infty ,p}))\Phi _{\infty
,p}(dz)\prod_{q=p+1}^{P}\Delta _{\sigma \xi _{\infty ,q}}(0,c_{q}\sigma \xi
_{\infty ,q})  \label{4.1}
\end{align}%
with $p_{\ast }=\max \{p_{0}(\theta ),\O \}$. Note that $G_{\infty ,\theta
,\sigma }(t)$ depends on $\theta $ only through $p_{\ast }$. Let $\hat{\Phi}%
_{n,p}$\ denote the cdf of a $k$-variate Gaussian random vector with mean
zero and variance-covariance matrix $\hat{\sigma}^{2}A[p](X[p]^{\prime
}X[p]/n)^{-1}A[p]^{\prime }$, $0<p\leq P$; we also adopt the convention that 
$\hat{\Phi}_{n,0}$ denotes the cdf of point-mass at zero in $\mathbf{R}^{k}$%
. [We use the same convention for $\hat{\Phi}_{n,p}$ in case $\hat{\sigma}=0$%
, which is a probability zero event.] An estimator $\check{G}_{n}(t)$ of $%
G_{n,\theta ,\sigma }(t)$\ is now defined as follows: We first employ an
auxiliary procedure $\bar{p}$ that consistently estimates $p_{0}(\theta )$
(e.g., $\bar{p}$ could be obtained from BIC or from a `general-to-specific'
hypothesis testing procedure employing critical values that go to infinity
but are $o(n^{1/2})$ as $n\rightarrow \infty $). The estimator $\check{G}%
_{n}(t)$ is now given by the expression in (\ref{4.1}) but with\ $p_{\ast }$%
, $\sigma $, $b_{\infty ,p}$, $\zeta _{\infty ,p}$, $\xi _{\infty ,p}$, and $%
\Phi _{\infty ,p}$ replaced by $\max \{\bar{p},\O \}$, $\hat{\sigma}$, $%
b_{n,p}$, $\zeta _{n,p}$, $\xi _{n,p}$, and $\hat{\Phi}_{n,p}$,
respectively. A little reflection shows that $\check{G}_{n}$ is again a cdf.
We have the following consistency results.

\begin{proposition}
\label{p4.1.1}The estimator $\check{G}_{n}$ is consistent (in the total
variation distance) for $G_{n,\theta ,\sigma }$ and $G_{\infty ,\theta
,\sigma }$. That is, for every $\delta >0$ 
\begin{equation}
P_{n,\theta ,\sigma }\left( \left\vert \left\vert \check{G}_{n}(\cdot
)-G_{n,\theta ,\sigma }(\cdot )\right\vert \right\vert _{TV}\,>\,\delta
\,\right) \;\overset{n\rightarrow \infty }{\longrightarrow }\;0,  \label{4.2}
\end{equation}%
\begin{equation}
P_{n,\theta ,\sigma }\left( \left\vert \left\vert \check{G}_{n}(\cdot
)-G_{\infty ,\theta ,\sigma }(\cdot )\right\vert \right\vert
_{TV}\,>\,\delta \,\right) \;\overset{n\rightarrow \infty }{\longrightarrow }%
\;0  \label{4.3}
\end{equation}%
for all $\theta \in \mathbf{R}^{P}$ and all $\sigma >0$.
\end{proposition}

While the estimator constructed above on the basis of the formula for $%
G_{\infty ,\theta ,\sigma }$ is consistent, it can be expected to perform
poorly in finite samples since convergence of $G_{n,\theta ,\sigma }$ to $%
G_{\infty ,\theta ,\sigma }$ is typically not uniform in $\theta $ (cf.
Appendix \ref{ae}), and since in case the true $\theta $ is `close' to $%
M_{p_{0}(\theta )-1}$ the auxiliary decision procedure $\bar{p}$ (although
being consistent for $p_{0}(\theta )$) will then have difficulties making
the correct decision in finite samples. In the next section we show that
this poor performance is not particular to the estimator $\check{G}_{n}$
constructed above, but is a genuine feature of the estimation problem under
consideration.

\subsubsection{Performance Limits and Impossibility Results\label{s5}}

We now provide lower bounds for the performance of estimators of the cdf $%
G_{n,\theta ,\sigma }(t)$\ of the post-model-selection estimator $A\tilde{%
\theta}$; that is, we give lower bounds on the worst-case probability that
the estimation error exceeds a certain threshold. These lower bounds are
large, being $1$ or $1/2$, depending on the situation considered;
furthermore, they remain lower bounds even if one restricts attention only
to certain subsets of the parameter space that shrink at the rate $n^{-1/2}$%
. In this sense the `impossibility' results are of a local nature. In
particular, the lower bounds imply that no uniformly consistent estimator of
the cdf $G_{n,\theta ,\sigma }(t)$\ exists, not even locally.

In the following, the asymptotic `correlation' between $A\tilde{\theta}(p)$
and $\tilde{\theta}_{p}(p)$ as measured by $C_{\infty
}^{(p)}=\lim_{n\rightarrow \infty }C_{n}^{(p)}$ will play an important r\^{o}%
le. [Recall that $\tilde{\theta}(p)$ denotes the least-squares estimator of $%
\theta $ based on model $M_{p}$ and that $A\theta $ is the parameter vector
of interest. Furthermore, the vector of covariances between $A\tilde{\theta}%
(p)$ and $\tilde{\theta}_{p}(p)$ is given by $\sigma ^{2}n^{-1}C_{n}^{(p)}$
with $C_{n}^{(p)}=A[p](X[p]^{\prime }X[p]/n)^{-1}e_{p}$.] Note that $%
C_{\infty }^{(p)}$ equals $A[p]Q[p:p]^{-1}e_{p}$, and hence does \emph{not}
depend on the unknown parameters $\theta $ or $\sigma $. In the important
special case discussed in the Introduction, cf. (\ref{1.1}), the matrix $A$
equals the $\O \times P$ matrix $(I_{\O }:0)$, and the condition $C_{\infty
}^{(p)}\neq 0$ reduces to the condition that the regressor corresponding to
the $p$-th column of $(V:W)$ is asymptotically correlated with at least one
of the regressors corresponding to the columns of $V$. See Example 1 below
for more discussion.

In the result to follow we shall consider performance limits for estimators
of $G_{n,\theta ,\sigma }(t)$ at a \emph{fixed} value of the argument $t$.
An estimator of $G_{n,\theta ,\sigma }(t)$ is now nothing else than a
real-valued random variable $\Gamma _{n}=\Gamma _{n}(Y,X)$. For mnemonic
reasons we shall, however, use the symbol $\hat{G}_{n}(t)$\ instead of $%
\Gamma _{n}$ to denote an arbitrary estimator of $G_{n,\theta ,\sigma }(t)$.
This notation should not be taken as implying that the estimator is obtained
by evaluating an estimated cdf at the argument $t$, or that it is a priori
constrained to lie between zero and one. We shall use this notational
convention mutatis mutandis also in subsequent sections. Regarding the
non-uniformity phenomenon, we then have a dichotomy which is described in
the following two results.

\begin{theorem}
\label{t5.1} Suppose that $A\tilde{\theta}(q)$ and $\tilde{\theta}_{q}(q)$
are asymptotically correlated, i.e., $C_{\infty }^{(q)}\neq 0$, for some $q$
satisfying $\O <{q}\leq P$, and let $q^{\ast }$ denote the largest $q$ with
this property. Then the following holds for every $\theta \in M_{q{^{\ast }}%
-1}$, every $\sigma $, $0<\sigma <\infty $, and every $t\in \mathbf{R}^{k}$:
There exist $\delta _{0}>0$ and $\rho _{0}$, $0<\rho _{0}<\infty $, such
that any estimator $\hat{G}_{n}(t)$ of $G_{n,\theta ,\sigma }(t)$ satisfying%
\begin{equation}
P_{n,\theta ,\sigma }\left( \left\vert \hat{G}_{n}(t)-G_{n,\theta ,\sigma
}(t)\right\vert >\delta \right) \overset{n\rightarrow \infty }{%
\longrightarrow }0  \label{5.1}
\end{equation}%
for each $\delta >0$ (in particular, every estimator that is consistent)
also satisfies 
\begin{equation}
\sup_{\overset{\vartheta \in M_{q^{\ast }}}{||\vartheta -\theta ||<\rho _{0}/%
\sqrt{n}}}P_{n,\vartheta ,\sigma }\left( \left\vert \hat{G}%
_{n}(t)-G_{n,\vartheta ,\sigma }(t)\right\vert >\delta _{0}\right) \overset{%
n\rightarrow \infty }{\longrightarrow }1.  \label{5.2}
\end{equation}%
The constants $\delta _{0}$ and $\rho _{0}$ may be chosen in such a way that
they depend only on $t$, $Q$, $A$, $\sigma $, and the critical values $c_{p}$
for $\O <{p}\leq P$. Moreover,%
\begin{equation}
\liminf_{n\rightarrow \infty }\inf_{\hat{G}_{n}(t)}\sup_{\overset{\vartheta
\in M_{q^{\ast }}}{||\vartheta -\theta ||<\rho _{0}/\sqrt{n}}}P_{n,\vartheta
,\sigma }\left( \left\vert \hat{G}_{n}(t)-G_{n,\vartheta ,\sigma
}(t)\right\vert >\delta _{0}\right) >0  \label{5.3}
\end{equation}%
and%
\begin{equation}
\sup_{\delta >0}\liminf_{n\rightarrow \infty }\inf_{\hat{G}_{n}(t)}\sup_{%
\overset{\vartheta \in M_{q^{\ast }}}{||\vartheta -\theta ||<\rho _{0}/\sqrt{%
n}}}P_{n,\vartheta ,\sigma }\left( \left\vert \hat{G}_{n}(t)-G_{n,\vartheta
,\sigma }(t)\right\vert >\delta \right) \geq \frac{1}{2},  \label{5.4}
\end{equation}%
where the infima in (\ref{5.3}) and (\ref{5.4}) extend over \emph{all }%
estimators $\hat{G}_{n}(t)$ of $G_{n,\theta ,\sigma }(t)$.
\end{theorem}

\begin{remark}
\normalfont\label{r5.1}Assume that the conditions of the preceding theorem
are satisfied. Suppose further that $p_{\odot }$, $\O \leq p_{\odot
}<q^{\ast }$, is such that either $p_{\odot }>0$ and some row of $A[p_{\odot
}]$ equals zero, or such that $p_{\odot }=0$. Then there exist $\delta
_{0}>0 $ and $0<\rho _{0}<\infty $ such that the left-hand side of (\ref{5.3}%
) is not less than $1/2$ for each $\theta \in M_{p_{\odot }}$.
\end{remark}

Theorem \ref{t5.1} a fortiori implies a corresponding `impossibility' result
for estimation of the function $G_{n,\theta ,\sigma }(\cdot )$ when the
estimation error is measured in the total variation distance or the
sup-norm; cf. also Section \ref{after}.

It remains to consider the -- quite exceptional -- case where the assumption
of Theorem~\ref{t5.1} is not satisfied, i.e., where $C_{\infty }^{(q)}=0$,
for all $q$ in the range $\O <{q}\leq P$. Under this `uncorrelatedness'\
condition it is indeed possible to construct an estimator of $G_{n,\theta
,\sigma }$ which is uniformly consistent: It is not difficult to see that
the asymptotic distribution of $G_{n,\theta ,\sigma }$ reduces to $\Phi
_{\infty ,P}$ under this `uncorrelatedness'\ condition. Furthermore, the
second half of Proposition~\ref{p.ae1} in Appendix \ref{ae} shows that then
the convergence of $G_{n,\theta ,\sigma }$ to its large-sample limit is
uniform w.r.t. $\theta $, suggesting $\hat{\Phi}_{n,P}$, an estimated
version of $\Phi _{\infty ,P}$, as an estimator for $G_{n,\theta ,\sigma }$.

\begin{proposition}
\label{p5.2} Suppose that $A\tilde{\theta}(q)$ and $\tilde{\theta}_{q}(q)$
are asymptotically uncorrelated, i.e., $C_{\infty }^{(q)}=0$, for all $q$
satisfying $\O <{q}\leq P$. Then%
\begin{equation}
\sup_{\theta \in \mathbf{R}^{P}}\,\sup_{\overset{\sigma \in \mathbf{R}}{%
\sigma _{\ast }\leq \sigma \leq \sigma ^{\ast }}}\,P_{n,\theta ,\sigma
}\left( \left\vert \left\vert \hat{\Phi}_{n,P}-G_{n,\theta ,\sigma
}\right\vert \right\vert _{TV}>\delta \right) \overset{n\rightarrow \infty }{%
\longrightarrow }0  \label{5.6}
\end{equation}%
holds for each $\delta >0$, and for any constants $\sigma _{\ast }$ and $%
\sigma ^{\ast }$ satisfying $0<\sigma _{\ast }\leq \sigma ^{\ast }<\infty $.
\end{proposition}

Inspection of the proof of Proposition~\ref{p5.2} shows that (\ref{5.6})
continues to hold if the estimator $\hat{\Phi}_{n,P}$ is replaced by any of
the estimators $\hat{\Phi}_{n,p}$ for $\O \leq p\leq P$. We also note that
in case $\O =0$ the assumption of Proposition~\ref{p5.2} is never satisfied
in view of Proposition~4.4 in Leeb and P\"{o}tscher (2006b), and hence
Theorem~\ref{t5.1} always applies in that case. Another consequence of
Proposition~4.4 in Leeb and P\"{o}tscher (2006b) is that -- under the
`uncorrelatedness'\ assumption of Proposition~\ref{p5.2} -- the restricted
least squares estimators $A\tilde{\theta}(q)$ for $q\geq \O $ perform
asymptotically as well as the unrestricted estimator $A\tilde{\theta}(P)$;
this clearly shows that the case covered by Proposition~\ref{p5.2} is highly
exceptional.

In summary we see that it is typically impossible to construct an estimator
of $G_{n,\theta ,\sigma }(t)$ which performs reasonably well even
asymptotically. Whenever Theorem~\ref{t5.1} applies, any estimator of $%
G_{n,\theta ,\sigma }(t)$ suffers from a non-uniformity defect which is
caused by parameters belonging to shrinking `tubes' surrounding $M_{q^{\ast
}-1}$. For the sake of completeness, we remark that outside a `tube' of
fixed positive radius that surrounds $M_{q^{\ast }-1}$ the non-uniformity
need not be present: Let $q^{\ast }$ be as in Theorem~\ref{t5.1} and define
the set $U$ as $U=\{\theta \in \mathbf{R}^{P}:\,|\theta _{q^{\ast }}|\geq
r\} $ for some fixed $r>0$. Then $\hat{\Phi}_{n,P}(t)$ is an estimator of $%
G_{n,\theta ,\sigma }(t)$ that is uniformly consistent over $\theta \in U$;
more generally, it can be shown that then the relation (\ref{5.6}) holds if
the supremum over $\theta $ on the left-hand side is restricted to $\theta
\in U$.

We conclude this section by illustrating the above results with some
important examples.

\textbf{Example 1: }\textit{(The distribution of }$\tilde{\chi}$) Consider
the model given in (\ref{1.1}) with $\chi $ representing the parameter of
interest. Using the general notation of Section \ref{s2}, this corresponds
to the case $A\theta =(\theta _{1},\dots ,\theta _{\O })^{\prime }=\chi $
with $A$ representing the $\O \times P$ matrix $(I_{\O }:0)$. Here $k=\O >0$%
. The cdf $G_{n,\theta ,\sigma }$ then represents the cdf of $\sqrt{n}\left( 
\tilde{\chi}-\chi \right) $. Assume first that $\lim_{n\rightarrow \infty
}V^{\prime }W/n\neq 0$. Then $C_{\infty }^{(q)}\neq 0$ holds for some $q>\O $%
. Consequently, the `impossibility' results for the estimation of $%
G_{n,\theta ,\sigma }$ given in Theorem~\ref{t5.1} always apply. Next assume
that $\lim_{n\rightarrow \infty }V^{\prime }W/n=0$. Then $C_{\infty
}^{(q)}=0 $ for every $q>\O $. In this case Proposition~\ref{p5.2} applies
and a uniformly consistent estimator of $G_{n,\theta ,\sigma }$ indeed
exists. Summarizing we note that any estimator of $G_{n,\theta ,\sigma }$
suffers from the non-uniformity phenomenon except in the special case where
the columns of $V$ and $W$ are asymptotically orthogonal in the sense that $%
\lim_{n\rightarrow \infty }V^{\prime }W/n=0$. But this is precisely the
situation where inclusion or exclusion of the regressors in $W$ has no
effect on the distribution of the estimator $\tilde{\chi}$ asymptotically;
hence it is not surprising that also the model selection procedure does not
have an effect on the estimation of the cdf of the post-model-selection
estimator $\tilde{\chi}$. This observation may tempt one to enforce
orthogonality between the columns of $V$ and $W$ by either replacing the
columns of $V$ by their residuals from the projection on the column space of 
$W$ or vice versa. However, this is not helpful for the following reasons:
In the first case one then in fact avoids model selection as all the
restricted least-squares estimators for $\chi $ under consideration (and
hence also the post-model selection estimator $\tilde{\chi}$) in the
reparameterized model coincide with the unrestricted least-squares
estimator. In the second case the coefficients of the columns of $V$ in the
reparameterized model no longer coincide with the parameter of interest $%
\chi $ (and again are estimated by one and the same estimator regardless of
inclusion/exclusion of columns of the transformed $W$-matrix).

\textbf{Example 2:} \textit{(The distribution of }$\tilde{\theta}$\textit{)}
For $A$ equal to $I_{P}$, the cdf $G_{n,\theta ,\sigma }$ is the cdf of $%
\sqrt{n}(\tilde{\theta}-\theta )$. Here, $A\tilde{\theta}(q)$ reduces to $%
\tilde{\theta}(q)$, and hence $A\tilde{\theta}(q)$ and $\tilde{\theta}%
_{q}(q) $ are perfectly correlated for every $q>\O $. Consequently, the
`impossibility' result for estimation of $G_{n,\theta ,\sigma }$ given in
Theorem~\ref{t5.1} applies. [In fact, the slightly stronger result mentioned
in Remark~\ref{r5.1}\ always applies here.] We therefore see that estimation
of the distribution of the post-model-selection estimator of the entire
parameter vector is always plagued by the non-uniformity phenomenon.

\textbf{Example 3:}\texttt{\ }\textit{(The distribution of a linear
predictor)} Suppose $A\neq 0$ is a $1\times P$ vector and one is interested
in estimating the cdf $G_{n,\theta ,\sigma }$ of the linear predictor $A%
\tilde{\theta}$. Then Theorem~\ref{t5.1} and the discussion following
Proposition~\ref{p5.2} show that the non-uniformity phenomenon always arises
in this estimation problem in case $\O =0$. In case $\O >0$, the
non-uniformity problem is generically also present, except in the degenerate
case where $C_{\infty }^{(q)}=0$, for all $q$ satisfying $\O <q\leq P$ (in
which case Proposition~4.4 in Leeb and P\"{o}tscher (2006b) shows that the
least-squares predictors from all models $M_{p}$, $\O \leq p\leq P$, perform
asymptotically equally well).

\section{Extensions to Other Model Selection Procedures Including AIC\label%
{s99}}

In this section we show that the `impossibility' result obtained in the
previous section for a `general-to-specific' model selection procedure
carries over to a large class of model selection procedures, including
Akaike's widely used AIC. Again consider the linear regression model (\ref%
{2.1}) with the same assumptions on the regressors and the errors as in
Section \ref{s2}. Let $\{0,1\}^{P}$ denote the set of all 0-1 sequences of
length $P$.$\mathfrak{\ }$For each $\mathfrak{r\in }\{0,1\}^{P}$ let $M_{%
\mathfrak{r}}$ denote the set $\{\theta \in \mathbf{R}^{P}:\,\theta _{i}(1-%
\mathfrak{r}_{i})=0$ $for$ $1\leq i\leq P\}$ where $\mathfrak{r}_{i}$
represents the $i$-th component of $\mathfrak{r}$. I.e., $M_{\mathfrak{r}}$
describes a linear submodel with those parameters $\theta _{i}$ restricted
to zero for which $\mathfrak{r}_{i}=0$. Now let $\mathfrak{R}$ be a
user-supplied subset of $\{0,1\}^{P}$. We consider model selection
procedures that select from the set $\mathfrak{R}$, or equivalently from the
set of models $\{M_{\mathfrak{r}}:\mathfrak{r\in R}\}$. Note that there is
now no assumption that the candidate models are nested (for example, if $%
\mathfrak{R}=\{0,1\}^{P}$ all possible submodels are candidates for
selection). Also cases where the inclusion of a subset of regressors is
undisputed on a priori grounds are obviously covered by this framework upon
suitable choice of $\mathfrak{R}$.

We shall assume throughout this section that $\mathfrak{R}$ contains $%
\mathfrak{r}_{full}=(1,\ldots ,1)$ and also at least one element $\mathfrak{r%
}_{\ast }$ satisfying $\left\vert \mathfrak{r}_{\ast }\right\vert =P-1$,
where $\left\vert \mathfrak{r}_{\ast }\right\vert $ represents the number of
non-zero coordinates of $\mathfrak{r}_{\ast }$. Let $\hat{\mathfrak{r}}$ be
an arbitrary model selection procedure, i.e., $\hat{\mathfrak{r}}=\hat{%
\mathfrak{r}}(Y,X)$ is a random variable taking its values in $\mathfrak{R}$%
. We furthermore assume throughout this section that the model selection
procedure $\hat{\mathfrak{r}}$ satisfies the following mild condition: For
every $\mathfrak{r}_{\ast }\in \mathfrak{R}$ with $\left\vert \mathfrak{r}%
_{\ast }\right\vert =P-1$ there exists a positive finite constant $c$
(possibly depending on $\mathfrak{r}_{\ast }$) such that for every $\theta
\in M_{\mathfrak{r}_{\ast }}$ which has exactly $P-1$ non-zero coordinates 
\begin{equation}
\lim_{n\rightarrow \infty }P_{n,\theta ,\sigma }\left( \{\hat{\mathfrak{r}}=%
\mathfrak{r}_{full}\}\blacktriangle \{\left\vert T_{\mathfrak{r}_{\ast
}}\right\vert \geq c\}\right) =\lim_{n\rightarrow \infty }P_{n,\theta
,\sigma }\left( \{\hat{\mathfrak{r}}=\mathfrak{r}_{\ast }\}\blacktriangle
\{\left\vert T_{\mathfrak{r}_{\ast }}\right\vert <c\}\right) =0  \label{99.1}
\end{equation}%
holds for every $0<\sigma <\infty $. Here $\blacktriangle $ denotes the
symmetric difference operator and $T_{\mathfrak{r}_{\ast }}$ represents the
usual t-statistic for testing the hypothesis $\theta _{i(\mathfrak{r}_{\ast
})}=0$ in the full model, where $i(\mathfrak{r}_{\ast })$ denotes the index
of the unique coordinate of $\mathfrak{r}_{\ast }$ that equals zero.

The above condition is quite natural for the following reason: For $\theta
\in M_{\mathfrak{r}_{\ast }}$ with exactly $P-1$ non-zero coordinates, every
reasonable model selection procedure will -- with probability approaching
unity -- decide only between $M_{\mathfrak{r}_{\ast }}$ and $M_{\mathfrak{r}%
_{full}}$; it is then quite natural that this decision will be based (at
least asymptotically) on the likelihood ratio between these two models,
which in turn boils down to the t-statistic. As will be shown below,
condition (\ref{99.1}) holds in particular for AIC-like procedures.

Let $A$ be a non-stochastic $k\times P$ matrix of full row rank $k$, $1\leq
k\leq P$, as in Section \ref{s3}. We then consider the cdf 
\begin{equation}
K_{n,\theta ,\sigma }(t)=P_{n,\theta ,\sigma }\left( \sqrt{n}A(\bar{\theta}%
-\theta )\leq t\right) \qquad (t\in \mathbf{R}^{k})  \label{99.2}
\end{equation}%
of a linear transformation of the post-model-selection estimator $\bar{\theta%
}$ obtained from the model selection procedure $\hat{\mathfrak{r}}$, i.e., 
\begin{equation*}
\bar{\theta}=\sum_{\mathfrak{r}\in \mathfrak{R}}\tilde{\theta}(\mathfrak{r}%
)\,\mathbf{1(}\hat{\mathfrak{r}}=\mathfrak{r})
\end{equation*}%
where the $P\times 1$\ vector $\tilde{\theta}(\mathfrak{r})$ represents the
restricted least-squares estimator obtained from model $M_{\mathfrak{r}}$,
with the convention that $\tilde{\theta}(\mathfrak{r})=0\in \mathbf{R}^{P}$
in case $\mathfrak{r}=(0,\ldots ,0)$. We then obtain the following result
for estimation of $K_{n,\theta ,\sigma }(t)$ at a \emph{fixed} value of the
argument $t$ which parallels the corresponding `impossibility' result in
Theorem \ref{t5.1}.

\begin{theorem}
\label{t99.2}Let $\mathfrak{r}_{\ast }\in \mathfrak{R}$ satisfy $\left\vert 
\mathfrak{r}_{\ast }\right\vert =P-1$, and let $i(\mathfrak{r}_{\ast })$
denote the index of the unique coordinate of $\mathfrak{r}_{\ast }$ that
equals zero; furthermore, let $c$ be the constant in (\ref{99.1})
corresponding to $\mathfrak{r}_{\ast }$. Suppose that $A\tilde{\theta}(%
\mathfrak{r}_{full})$ and $\tilde{\theta}_{i(\mathfrak{r}_{\ast })}(%
\mathfrak{r}_{full})$ are asymptotically correlated, i.e., $AQ^{-1}e_{_{i(%
\mathfrak{r}_{\ast })}}\neq 0$, where $e_{_{i(\mathfrak{r}_{\ast })}}$
denotes the $i(\mathfrak{r}_{\ast })$-th standard basis vector in $\mathbf{R}%
^{P}$. Then for every $\theta \in M_{\mathfrak{r}_{\ast }}$ which has
exactly $P-1$ non-zero coordinates, for every $\sigma $, $0<\sigma <\infty $%
, and for every $t\in \mathbf{R}^{k}$ the following holds: There exist $%
\delta _{0}>0$ and $\rho _{0}$, $0<\rho _{0}<\infty $, such that any
estimator $\hat{K}_{n}(t)$ of $K_{n,\theta ,\sigma }(t)$ satisfying 
\begin{equation}
P_{n,\theta ,\sigma }\left( \left\vert \hat{K}_{n}(t)-K_{n,\theta ,\sigma
}(t)\right\vert \,>\,\delta \right) \overset{n\rightarrow \infty }{%
\longrightarrow }0  \label{99.8}
\end{equation}%
for each $\delta >0$ (in particular, every estimator that is consistent)
also satisfies 
\begin{equation}
\sup_{\overset{\vartheta \in \mathbf{R}^{P}}{||\vartheta -\theta ||<\rho
_{0}/\sqrt{n}}}P_{n,\vartheta ,\sigma }\left( \left\vert \hat{K}%
_{n}(t)-K_{n,\vartheta ,\sigma }(t)\right\vert \,>\,\delta _{0}\right) 
\overset{n\rightarrow \infty }{\longrightarrow }1\,.  \label{99.9}
\end{equation}%
The constants $\delta _{0}$ and $\rho _{0}$ may be chosen in such a way that
they depend only on $t,Q,A,\sigma ,$ and $c$. Moreover, 
\begin{equation}
\liminf_{n\rightarrow \infty }\inf_{\hat{K}_{n}(t)}\sup_{\overset{\vartheta
\in \mathbf{R}^{P}}{||\vartheta -\theta ||<\rho _{0}/\sqrt{n}}%
}\,P_{n,\vartheta ,\sigma }\left( \left\vert \hat{K}_{n}(t)-K_{n,\vartheta
,\sigma }(t)\right\vert \,>\,\delta _{0}\right) >0  \label{99.10}
\end{equation}%
and 
\begin{equation}
\sup_{\delta >0}\,\liminf_{n\rightarrow \infty }\inf_{\hat{K}_{n}(t)}\sup_{%
\overset{\vartheta \in \mathbf{R}^{P}}{||\vartheta -\theta ||<\rho _{0}/%
\sqrt{n}}}\,P_{n,\vartheta ,\sigma }\left( \left\vert \hat{K}%
_{n}(t)-K_{n,\vartheta ,\sigma }(t)\right\vert \,>\,\delta \right) \,\geq 1/2
\label{99.11}
\end{equation}%
hold, where the infima in (\ref{99.10}) and (\ref{99.11}) extend over \emph{%
all }estimators $\hat{K}_{n}(t)$ of $K_{n,\theta ,\sigma }(t)$.
\end{theorem}

The basic condition (\ref{99.1}) on the model selection procedure employed
in the above result will certainly hold for any hypothesis testing procedure
that (i) asymptotically selects only correct models, (ii) employs a
likelihood ratio test (or an asymptotically equivalent test) for testing $M_{%
\mathfrak{r}_{full}}$ versus smaller models (at least versus the models $M_{%
\mathfrak{r}_{\ast }}$ with $\mathfrak{r}_{\ast }$ as in condition (\ref%
{99.1})), and (iii) uses a critical value for the likelihood ratio test that
converges to a finite positive constant. In particular, this applies to
usual thresholding procedures as well as to a variant of the
`general-to-specific' procedure discussed in Section \ref{s2} where the
error variance in the construction of the test statistic for hypothesis $%
H_{0}^{p}$ is estimated from the fitted model $M_{p}$ rather than from the
overall model. We next verify condition (\ref{99.1}) for AIC-like
procedures. Let $RSS(\mathfrak{r})$ denote the residual sum of squares from
the regression employing model $M_{\mathfrak{r}}$ and set%
\begin{equation}
IC(\mathfrak{r})=\log \left( RSS(\mathfrak{r})\right) +\left\vert \mathfrak{r%
}\right\vert \Upsilon _{n}/n  \label{99.12}
\end{equation}%
where $\Upsilon _{n}\geq 0$ denotes a sequence of real numbers satisfying $%
\lim_{n\rightarrow \infty }\Upsilon _{n}=\Upsilon $ and $\Upsilon $ is a
positive real number. Of course, $IC(\mathfrak{r})=AIC(\mathfrak{r})$ if $%
\Upsilon _{n}=2$. The model selection procedure $\hat{\mathfrak{r}}_{IC}$ is
then defined as a minimizer (more precisely, as a measurable selection from
the set of minimizers) of $IC(\mathfrak{r})$ over $\mathfrak{R}$. It is
well-known that the probability that $\hat{\mathfrak{r}}_{IC}$ selects an
incorrect model converges to zero. Hence, elementary calculations show that
condition (\ref{99.1}) is satisfied for $c=\Upsilon ^{1/2}$.

The analysis of post-model-selection estimators based on AIC-like model
selection procedures given in this section proceeded by bringing this case
under the umbrella of the results obtained in Section \ref{s2}. Verification
of condition (\ref{99.1}) is the key that enables this approach. A complete
analysis of post-model-selection estimators based on AIC-like model
selection procedures, similar to the analysis in Section \ref{s2} for the
`general-to-specific' model selection procedure, is certainly possible but
requires a direct and detailed analysis of the distribution of this
post-model-selection estimator. [Even the mild condition that $\mathfrak{R}$
contains $\mathfrak{r}_{full}$ and also at least one element $\mathfrak{r}%
_{\ast }$ satisfying $\left\vert \mathfrak{r}_{\ast }\right\vert =P-1$ can
then be relaxed in such an analysis.] We furthermore note that in the
special case where $\mathfrak{R}=\{\mathfrak{r}_{full},\mathfrak{r}_{\ast
}\} $ and an AIC-like model selection procedure as in (\ref{99.12}) is used,
the results in the above theorem in fact hold for all $\theta \in M_{%
\mathfrak{r}_{\ast }}$.

\section{Remarks and Extensions\label{s6}}

\begin{remark}
\normalfont\label{r6.0} Although not emphasized in the notation, all results
in the paper also hold if the elements of the design matrix $X$ depend on
sample size. Furthermore, all results are expressed solely in terms of the
distributions $P_{n,\theta ,\sigma }(\cdot )$ of $Y$, and hence they also
apply if the elements of $Y$ depend on sample size, including the case where
the random vectors $Y$ are defined on different probability spaces for
different sample sizes.
\end{remark}

\begin{remark}
\normalfont\label{r6.4} The model selection procedure considered in Section %
\ref{s2}\ is based on a sequence of tests which use critical values $c_{p}$
that do not depend on sample size and satisfy $0<c_{p}<\infty $ for $\O %
<p\leq P$. If these critical values are allowed to depend on sample size
such that they now satisfy $c_{n,p}\rightarrow c_{\infty ,p}$ as $%
n\rightarrow \infty $ with $0<c_{\infty ,p}<\infty $ for $\O <p\leq P$, the
results in Leeb and P\"{o}tscher (2003) as well as in Leeb (2005, 2006)
continue to hold; see Remark~6.2(i) in Leeb and P\"{o}tscher (2003) and
Remark~6.1(ii) in Leeb (2005). As a consequence, the results in the present
paper can also be extended to this case quite easily.
\end{remark}

\begin{remark}
\normalfont\label{r6.1} The `impossibility' results given in Theorems~\ref%
{t5.1} and \ref{t99.2} (as well as the variants thereof discussed in the
subsequent Remarks~\ref{r6.2}-\ref{r6.8}) also hold for the class of all
randomized estimators (with $P_{n,\theta ,\sigma }^{\ast }$ replacing $%
P_{n,\theta ,\sigma }$ in those results, where $P_{n,\theta ,\sigma }^{\ast
} $ denotes the distribution of the randomized sample). This follows
immediately from Lemma~3.6 and the attending discussion in Leeb and P\"{o}%
tscher (2006a).
\end{remark}

\begin{remark}
\normalfont\label{r6.2}

\begin{enumerate}
\item Let $\psi _{n,\theta ,\sigma }$ denote the expectation of $\tilde{%
\theta}$ under $P_{n,\theta ,\sigma }$, and consider the cdf $H_{n,\theta
,\sigma }(t)=P_{n,\theta ,\sigma }(\sqrt{n}A(\tilde{\theta}-\psi _{n,\theta
,\sigma })\leq t)$. Results for the cdf $H_{n,\theta ,\sigma }$\ quite
similar to the results for $G_{n,\theta ,\sigma }$\ obtained in the present
paper can be established. A similar remark applies to the
post-model-selection estimator $\bar{\theta}$ considered in Section \ref{s99}%
.

\item In Leeb (2006) also the cdf $G_{n,\theta ,\sigma }^{\ast }$ is
analyzed, which correspond to a (typically infeasible) model selection
procedure that makes use of knowledge of $\sigma $. Results completely
analogous to the ones in the present paper can also be obtained for this cdf.
\end{enumerate}
\end{remark}

\begin{remark}
\normalfont\label{r6.3} Results similar to the ones in Section \ref{s5} can
also be obtained for estimation of the asymptotic cdf $G_{\infty ,\theta
,\sigma }(t)$ (or of the asymptotic cdfs corresponding to the variants
discussed in the previous remark). Since these results are of limited
interest, we omit them. In particular, note that an `impossibility' result
for estimation of $G_{\infty ,\theta ,\sigma }(t)$ per se does \emph{not}
imply a corresponding `impossibility' result for estimation of $G_{n,\theta
,\sigma }(t)$, since $G_{n,\theta ,\sigma }(t)$ does in general not converge
uniformly to $G_{\infty ,\theta ,\sigma }(t)$ over the relevant subsets in
the parameter space; cf. Appendix \ref{ae}. [An analogous remark applies to
the model selection procedures considered in Section \ref{s99}.]
\end{remark}

\begin{remark}
\normalfont\label{r6.6} Let $\pi _{n,\theta ,\sigma }(p)$ denote the model
selection probability $P_{n,\theta ,\sigma }(\hat{p}=p)$, $\O \leq p\leq P$
corresponding to the model selection procedure discussed in Section \ref{s2}%
. The finite-sample properties and the large-sample limit behavior of these
quantities are thoroughly analyzed in Leeb (2006); cf. also Leeb and P\"{o}%
tscher (2003). For these model selection probabilities the following results
can be established which we discuss here only briefly:

\begin{enumerate}
\item The model selection probabilities $\pi _{n,\theta ,\sigma }(p)$
converge to well-defined large-sample limits which we denote by $\pi
_{\infty ,\theta ,\sigma }(p)$. Similar as in Proposition~\ref{p.ae1} in
Appendix \ref{ae}, the convergence of $\pi _{n,\theta ,\sigma }(p)$ to $\pi
_{\infty ,\theta ,\sigma }(p)$ is non-uniform w.r.t. $\theta $. [For the
case $\O =0$, this phenomenon is described in Corollary~5.6 of Leeb and P%
\"{o}tscher (2003).]

\item The model selection probabilities $\pi _{n,\theta ,\sigma }(p)$\ can
be estimated consistently. However, uniformly consistent estimation is again
not possible. A similar remark applies to the large-sample limits $\pi
_{\infty ,\theta ,\sigma }(p)$.
\end{enumerate}
\end{remark}

\begin{remark}
\normalfont\label{r6.8} `Impossibility' results similar to the ones given in
Theorems \ref{t5.1} and \ref{t99.2} for the cdf can also be obtained for
other characteristics of the distribution of a linear function of a
post-model-selection estimator like the mean-squared error or the bias of $%
\sqrt{n}A\tilde{\theta}$.
\end{remark}

\section{On the Scope of the Impossibility Results\label{after}}

The non-uniformity phenomenon described, e.g., in (\ref{5.2}) of Theorem \ref%
{t5.1} is caused by a mechanism that can informally be described as follows.
Under the assumptions of that theorem, one can find an appropriate $\theta $
and an appropriate sequence $\vartheta _{n}=\theta +\gamma /\sqrt{n}$
exhibiting two crucial properties:

\begin{enumerate}
\item \label{pa} The probability measures $P_{n,\vartheta _{n},\sigma }$
corresponding to $\vartheta _{n}$ are `close' to the measures $P_{n,\theta
,\sigma }$ corresponding to $\theta $, in the sense of contiguity. This
entails that an estimator, that converges to some limit in probability under 
$P_{n,\theta ,\sigma }$, converges to the \emph{same limit} also under $%
P_{n,\vartheta _{n},\sigma }$.

\item \label{pb} For given $t$, the estimands $G_{n,\vartheta _{n},\sigma
}(t)$ corresponding to $\vartheta _{n}$ are `far away' from the estimands $%
G_{n,\theta ,\sigma }(t)$ corresponding to $\theta $, in the sense that $%
G_{n,\vartheta _{n},\sigma }(t)$ and $G_{n,\theta ,\sigma }(t)$ converge to 
\emph{different limits}, i.e., $G_{\infty ,\theta ,\sigma ,0}(t)$ is
different from $G_{\infty ,\theta ,\sigma ,\gamma }(t)$.
\end{enumerate}

In view of Property \ref{pa}, an estimator $\hat{G}_{n}(t)$ satisfying $\hat{%
G}_{n}(t)-G_{n,\theta ,\sigma }(t)\rightarrow 0$ in probability under $%
P_{n,\theta ,\sigma }$, also satisfies $\hat{G}_{n}(t)-G_{n,\theta ,\sigma
}(t)\rightarrow 0$ in probability under $P_{n,\vartheta _{n},\sigma }$. In
view of Property \ref{pb}, such an estimator $\hat{G}_{n}(t)$ is hence `far
away' from the estimand $G_{n,\vartheta _{n},\sigma }(t)$ with high
probability under $P_{n,\vartheta _{n},\sigma }$. In other words, an
estimator that is close to $G_{n,\theta ,\sigma }(t)$ under $P_{n,\theta
,\sigma }$ must be far away from $G_{n,\vartheta _{n},\sigma }(t)$ under $%
P_{n,\vartheta _{n},\sigma }$. Formalized and refined, this argument leads
to (\ref{5.2}) and, as a consequence, to the non-existence of uniformly
consistent estimators for $G_{n,\theta ,\sigma }(t)$. [There are a number of
technical details in this formalization process that need careful attention
in order to obtain the results in their full strength as given in Sections %
\ref{s2} and \ref{s99}.]

The above informal argument that derives (\ref{5.2}) from Properties \ref{pa}
and \ref{pb} can be refined and formalized in a much more general and
abstract framework, see Section 3 of Leeb and P\"{o}tscher (2006a) and the
references therein. That paper also provides a general framework for
deriving results like (\ref{5.3}) and (\ref{5.4}) of Theorem \ref{t5.1}. The
mechanism leading to such lower bounds is similar to the one outlined above,
where for some of the results the concept of contiguity of the probability
measures involved has to be replaced by closeness of these measures in total
variation distance. We use the results in Section 3 of Leeb and P\"{o}tscher
(2006a) to formally convert Properties \ref{pa} and \ref{pb} into the
`impossibility' results of the present paper; cf. Appendix \ref{ac}.

Verifying the aforementioned Property \ref{pa} in the context of the present
paper is straightforward because we consider a Gaussian linear model. What
is technically more challenging and requires some work is the verification
of Property \ref{pb}; this is done in Appendix \ref{aa} inter alia and rests
on results of Leeb (2002, 2005, 2006).

Two important observations on Properties \ref{pa} and \ref{pb} are in order:
First, Property \ref{pa} is typically satisfied in general parametric models
under standard regularity conditions; e.g., it is satisfied whenever the
model is locally asymptotically normal. Second, Property \ref{pb} relies on
limiting properties only and not on the finite-sample structure of the
underlying statistical model. Now, the limit distributions of
post-model-selection estimators in sufficiently regular parametric or
semi-parametric models are typically the same as the limiting distributions
of the corresponding post-model-selection estimators in a Gaussian linear
model (see, e.g., Sen (1979), P\"{o}tscher (1991), Nickl (2003), or Hjort
and Claeskens (2003)). Hence, establishing Property \ref{pb} for the
Gaussian linear model then typically establishes the same result for a large
class of general parametric or semi-parametric models.\footnote{%
Some care has to be taken here. In the Gaussian linear case the
finite-sample cdfs converge at every value of the argument $t$, cf.
Propisition \ref{p3.1}. In a general parametric model, sometimes the
asymptotic results (e.g., Hjort and Claeskens (2003, Theorem 4.1)) only
guarantee weak convergence. Hence, to ensure convergence of the relevant
cdfs at a given argument $t$ as required in Proberty \ref{pb}, additional
considerations have to be employed. [This is, however, of no concern in the
context discussed in the next but one paragraph in this section.]
\par
{}} For example, Property \ref{pb} can be verified for a large class of
pre-test estimators in sufficiently regular parametric models by arguing as
in Appendix \ref{aa} and using the results of Nickl (2003) to reduce to the
Gaussian linear case. Hence, the impossibility result given in Theorem \ref%
{t5.1} can be extended to more general parametric and semiparametric models
with ease. The fact that we use a Gaussian linear model for the analysis in
the present paper is a matter of convenience rather than a necessity.

The non-uniformity results in Theorem \ref{t5.1} are for (conservative)
`general-to-specific' model selection from a nested family of models.
Theorem \ref{t99.2} extends this to more general (conservative) model
selection procedures (including AIC and related procedures) and to more
general families of models. The proof of Theorem \ref{t99.2} proceeds by
reducing the problem to one where only two nested models are considered, and
then to appeal to the results of Theorem \ref{t5.1}. The condition on the
model selection procedures that enables this reduction is condition (\ref%
{99.1}). It is apparent from the discussion in Section \ref{s99} that this
condition is satisfied for many model selection procedures. Furthermore, for
the same reasons as given in the preceding paragraph, also Theorem \ref%
{t99.2} can easily be extended to sufficiently regular parametric and
semi-parametric models.

The `impossibility' results in the present paper are formulated for
estimating $G_{n,\theta ,\sigma }(t)$ for a given value of $t$. Suppose that
we are now asking the question whether the cdf $G_{n,\theta ,\sigma }(\cdot
) $ viewed as a \emph{function} can be estimated uniformly consistently,
where consistency is relative to a metric that metrizes weak convergence.%
\footnote{%
Or, in fact, any metric w.r.t. which the relevant cdfs converge.} Using a
similar reasoning as above (which can again be made formal by using, e.g.,
Lemma 3.1 in Leeb and P\"{o}tscher (2006a)) the key step now is to show that
the\emph{\ function} $G_{\infty ,\theta ,\sigma ,0}(\cdot )$ is different
from the \emph{function} $G_{\infty ,\theta ,\sigma ,\gamma }(\cdot )$.
Obviously, it is a much simpler problem to find a $\gamma $ such that the 
\emph{functions }$G_{\infty ,\theta ,\sigma ,0}(\cdot )$ and $G_{\infty
,\theta ,\sigma ,\gamma }(\cdot )$ differ, than to find a $\gamma $ such
that the \emph{values} $G_{\infty ,\theta ,\sigma ,0}(t)$ and $G_{\infty
,\theta ,\sigma ,\gamma }(t)$ for a given $t$ differ. Certainly, having
solved the latter problem in Appendix \ref{aa}, this also provides an answer
to the former. This then immediately delivers the desired `impossibility'
result. [We note that in some special cases simpler arguments than the ones
used in Appendix \ref{aa} can be employed to solve the former problem: For
example, in case $A=I$ the functions $G_{\infty ,\theta ,\sigma ,0}(\cdot )$
and $G_{\infty ,\theta ,\sigma ,\gamma }(\cdot )$ can each be shown to be
convex combinations of cdfs that are concentrated on subspaces of different
dimensions. This can be exploited to establish without much difficulty that
the \emph{functions }$G_{\infty ,\theta ,\sigma ,0}(\cdot )$ and $G_{\infty
,\theta ,\sigma ,\gamma }(\cdot )$ differ. For purpose of comparison we note
that for general $A$ the distributions $G_{\infty ,\theta ,\sigma ,0}$ and $%
G_{\infty ,\theta ,\sigma ,\gamma }$ can both be absolutely continuous
w.r.t. Lebesgue measure, not allowing one to use this simple argument.]
Again the discussion in this paragraph extends to more general parametric
and semiparametric models without difficulty.

The present paper, including the discussion in this section, has focussed on
conservative model selection procedures. However, the discussion should make
it clear that similar `impossibility' results plague consistent model
selection. Section 2.3 in Leeb and P\"{o}tscher (2006a) in fact gives such
an `impossibility' result in a simple case.

We close with the following observations. Verification of Property \ref{pb},
whether it is for $G_{\infty ,\theta ,\sigma ,0}(t)$ and $G_{\infty ,\theta
,\sigma ,\gamma }(t)$ (for given $t)$ or for $G_{\infty ,\theta ,\sigma
,0}(\cdot )$ and $G_{\infty ,\theta ,\sigma ,\gamma }(\cdot )$, shows that
the post-model-selection estimator $A\tilde{\theta}$ is a so-called
non-regular estimator for $A\theta $: Consider an estimator $\tilde{\beta}$
in a parametric model $\{P_{n,\beta }:\beta \in B\}$ where the parameter
space $B$ is an open subset of Euclidean space $\mathbb{R}^{d}$. Suppose $%
\tilde{\beta}$, properly scaled and centered, has a limit distribution under
local alternatives, in the sense that $\sqrt{n}(\tilde{\beta}-(\beta +\gamma
/\sqrt{n}))$ converges in law under $P_{n,\beta +\gamma /\sqrt{n}}$ to a
limit distribution $L_{\infty ,\beta ,\gamma }(\cdot )$ for every $\gamma $.
The estimator $\tilde{\beta}$ is called \emph{regular} if for every $\beta $
the limit distribution $L_{\infty ,\beta ,\gamma }(\cdot )$ does not depend
on $\gamma $; cf. van der Vaart (1998, Section 8.5). Suppose now that the
model is, e.g., locally asymptotically normal (hence the contiguity property
in Property \ref{pa} is satisfied). The informal argument outlined at the
beginning of this section (and which is formalized in Lemma 3.1 of Leeb and P%
\"{o}tscher (2006a)) then in fact shows that the cdf of \emph{any}
non-regular estimator can not be estimated uniformly consistently (where
consistency is relative to any metric that metrizes weak convergence).

\section{Conclusions\label{s7}}

Despite the fact that we have shown that consistent estimators for the
distribution of a post-model-selection estimator can be constructed with
relative ease, we have also demonstrated that no estimator of this
distribution can have satisfactory performance (locally) uniformly in the
parameter space, even asymptotically. In particular, no (locally) uniformly
consistent estimator of this distribution exists. Hence, the answer to the
question posed in the title has to be negative. The results in the present
paper also cover the case of linear functions (e.g., predictors) of the
post-model-selection estimator.

We would like to stress here that resampling procedures like, e.g., the
bootstrap or subsampling, do not solve the problem at all. First note that
standard bootstrap techniques will typically not even provide consistent
estimators of the finite-sample distribution of the post-model-selection
estimator, as the bootstrap can be shown to stay random in the limit
(Kulperger and Ahmed (1992), Knight (1999, Example 3))\footnote{%
Brownstone (1990) claims the validity of a bootstrap procedure that is based
on a conservative model selection procedure in a linear regression model.
Kilian (1998) makes a similar claim in the context of autoregressive models
selected by a conservative model selection procedure. Also Hansen (2003)
contains such a claim for a stationary bootstrap procedure based on a
conservative model selection procedure. The above discussion intimates that
these claims are at least unsubstantiated.}. Basically the only way one can
coerce the bootstrap into delivering a consistent estimator is to resample
from a model that has been selected by an \emph{auxiliary} consistent model
selection procedure. The consistent estimator constructed in Section \ref%
{s4.1} is in fact of this form. In contrast to the standard bootstrap,
subsampling will typically deliver consistent estimators. However, the
`impossibility' results given in this paper apply to \emph{any }estimator
(including randomized estimators) of the cdf of a post-model-selection
estimator. Hence, also \emph{any} resampling based estimator suffers from
the non-uniformity defects described in Theorems \ref{t5.1} and \ref{t99.2};
cf. also Remark \ref{r6.1}.

The `impossibility' results in Theorems \ref{t5.1} and \ref{t99.2} are
derived in the framework of a normal linear regression model (and a fortiori
these results continue to hold in any model which includes the normal linear
regression model as a special case), but this is more a matter of
convenience than anything else: As discussed in Section \ref{after}, similar
results can be obtained in general statistical models allowing for
nonlinearity or dependent data, e.g., as long as standard regularity
conditions for maximum likelihood theory are satisfied.

The results in the present paper are derived for a large class of
conservative model selection procedures (i.e., procedures that select
overparameterized models with positive probability asymptotically) including
Akaike's AIC and typical `general-to-specific' hypothesis testing
procedures. For consistent model selection procedures -- like BIC or testing
procedures with suitably diverging critical values $c_{p}$ (cf. Bauer, P\"{o}%
tscher, and Hackl (1988)) -- the (pointwise) asymptotic distribution is
always normal. [This is elementary, cf. Lemma~1 in P\"{o}tscher (1991).]
However, as discussed at length in Leeb and P\"{o}tscher (2005a), this
asymptotic normality result paints a misleading picture of the finite sample
distribution which can be far from a normal, the convergence of the
finite-sample distribution to the asymptotic normal distribution not being
uniform. `Impossibility' results similar to the ones presented here can also
be obtained for post-model-selection estimators based on consistent model
selection procedures. These will be discussed in detail elsewhere. For a
simple special case such an `impossibility' result is given in Section 2.3
of Leeb and P\"{o}tscher (2006a).

The `impossibility' of estimating the distribution of the
post-model-selection estimator does not per se preclude the possibility of
conducting valid inference after model selection, a topic that deserves
further study. However, it certainly makes this a more challenging task.

\setlength{\parskip}{2pt} \appendix

\section{Auxiliary Lemmas \label{aa}}

\begin{lemma}
\label{la.1} Let $Z$ be a random vector with values in $\mathbf{R}^{k}$ and
let $W$ be a univariate standard Gaussian random variable independent of $Z$%
. Furthermore, let $C\in \mathbf{R}^{k}$ and $\tau >0$. Then 
\begin{equation}
\mathbb{P}(Z\leq Cx)\mathbb{P}(|W-x|<\tau )\,+\,\mathbb{P}(Z\leq
CW,|W-x|\geq \tau )  \label{la.1-1}
\end{equation}%
is constant as a function of $x\in \mathbf{R}$ if and only if $C=0$ or $%
\mathbb{P}(Z\leq Cx)=0$ for each $x\in \mathbf{R}$.
\end{lemma}

\textbf{Proof of Lemma \ref{la.1}:} Suppose $C=0$ holds. Using independence
of $Z$ and $W$ it is then easy to see that (\ref{la.1-1}) reduces to $%
\mathbb{P}(Z\leq 0)$, which is constant in $x$. If $\mathbb{P}(Z\leq Cx)=0$
for every $x\in \mathbf{R}$, then $\mathbb{P}(Z\leq CW)=0$, and hence (\ref%
{la.1-1}) is again constant, namely equal to zero.

To prove the converse, assume that (\ref{la.1-1}) is constant in $x\in 
\mathbf{R}$. Letting $x\rightarrow \infty $, we see that (\ref{la.1-1}) must
be equal to $\mathbb{P}(Z\leq CW)$. This entails that 
\begin{equation}
\mathbb{P}(Z\leq Cx)\mathbb{P}(|W-x|<\tau )=\mathbb{P}(Z\leq CW,\,|W-x|<\tau
)  \notag
\end{equation}%
holds for every $x\in \mathbf{R}$. Write $F(x)$ as shorthand for $\mathbb{P}%
(Z\leq Cx)$, and let $\Phi (z)$ and $\phi (z)$ denote the cdf and density of 
$W$, respectively. Then the expression in the above display can be written
as 
\begin{equation}
F(x)(\Phi (x+\tau )-\Phi (x-\tau ))=\int_{x-\tau }^{x+\tau }F(z)\phi
(z)dz.\qquad (x\in \mathbf{R})  \label{pla.1-2}
\end{equation}%
We now further assume that $C\neq 0$ and that $F(x)\neq 0$ for at least one $%
x\in \mathbf{R}$, and show that this leads to a contradiction.

Consider first the case where all components of $C$ are non-negative. Since $%
F$ is not identically zero, it is then, up to a scale factor, the cdf of a
random variable on the real line. But then (\ref{pla.1-2}) can not hold for
all $x\in \mathbf{R}$ as shown in Example 7 in Leeb (2002) (cf. also
equation (7) in that paper). The case where all components of $C$ are
non-positive follows similarly by applying the above argument to $F(-x)$ and
upon observing that both $\Phi (x+\tau )-\Phi (x-\tau )$ and $\phi (x)$ are
symmetric around $x=0$.

Finally, consider the case where $C$ has at least one positive and one
negative component. In this case clearly $\lim_{x\rightarrow -\infty
}F(x)=\lim_{x\rightarrow \infty }F(x)=0$ holds. Since $F(x)$ is continuous
in view of (\ref{pla.1-2}), we see that $F(x)$ attains its (positive)
maximum at some point $x_{1}\in \mathbf{R}$. Now note that (\ref{pla.1-2})
with $x_{1}$ replacing $x$ can be written as 
\begin{equation*}
\int_{x_{1}-\tau }^{x_{1}+\tau }(F(x_{1})-F(z))\phi (z)dz=0.
\end{equation*}%
This immediately entails that $F(x)=F(x_{1})$ for each $x\in \lbrack
x_{1}-\tau ,x_{1}+\tau ]$ (because $F(x)$ is continuous and because of the
definition of $x_{1}$). Repeating this argument with $x_{1}-\tau $ replacing 
$x_{1}$ and proceeding inductively, we obtain that $F(x)=F(x_{1})$ for each $%
x$ satisfying $x\leq x_{1}+\tau $, a contradiction with $\lim_{x\rightarrow
-\infty }F(x)=0$. {\hfill $\Box $}

\begin{lemma}
\label{la.2} Let $M$ and $N$ be matrices of dimension $k\times p$ and $%
k\times q$, respectively, such that the matrix $(M:N)$ has rank $k$ ($k\geq
1 $, $p\geq 1$, $q\geq 1$). Let $t\in \mathbf{R}^{k}$, and let $V$ be a
random vector with values in $\mathbf{R}^{p}$ whose distribution assigns
positive mass to every (non-empty) open subset of $\mathbf{R}^{p}$ (e.g., it
possesses an almost everywhere positive Lebesgue density). Set $f(x)=\mathbb{%
P}(MV\leq t+Nx)$, $x\in \mathbf{R}^{q}$. If one of the rows of $M$ consists
of zeros only, then $f$ is discontinuous at some point $x_{0}$. More
precisely, there exist $x_{0}\in \mathbf{R}^{q}$, $z\in \mathbf{R}^{q}$ and
a constant $c>0$, such that $f(x_{0}+\delta z)\geq c$ and $f(x_{0}-\delta
z)=0$ hold for every sufficiently small $\delta >0$.
\end{lemma}

\textbf{Proof of Lemma \ref{la.2}:} The case where $M$ is the zero-matrix is
trivial. Otherwise, let $I_{0}$ denote the set of indices $i$, $1\leq i\leq
k $, for which the $i$-th row of $M$ is zero. Let $(M_{0}:N_{0})$ denote the
matrix consisting of those rows of $(M:N)$ whose index is in $I_{0}$, and
let $(M_{1}:N_{1})$ denote the matrix consisting of the remaining rows of $%
(M:N)$. Clearly, $M_{0}$ is then the zero matrix. Furthermore, note that $%
N_{0}$ has full row-rank. Moreover, let $t_{0}$ denote the vector consisting
of those components of $t$ whose index is in $I_{0}$ and let $t_{1}$ denote
the vector containing the remaining components. With this notation, $f(x)$
can be written as $\mathbb{P}(0\leq t_{0}+N_{0}x,\;M_{1}V\leq t_{1}+N_{1}x)$.

For vectors $\mu \in \mathbf{R}^{p}$ and $\eta \in \mathbf{R}^{q}$ to be
specified in a moment, set $t^{\ast }=t+M\mu +N\eta $, and let $t_{0}^{\ast
} $ and $t_{1}^{\ast }$ be defined similarly to $t_{0}$ and $t_{1}$. Since
the matrix $(M:N)$ has full rank $k$, we can choose $\mu $ and $\eta $ such
that $t_{0}^{\ast }=0$ and $t_{1}^{\ast }>0$. Choose $z\in \mathbf{R}^{q}$
such that $N_{0}z>0$, which is possible because $N_{0}$ has full row-rank.
Set $x_{0}=\eta $. Then for every $\epsilon \in \mathbb{R}$ we have%
\begin{eqnarray*}
f(x_{0}+\epsilon z) &=&f(\eta +\epsilon z)=\mathbb{P}(MV\leq t+N(\eta
+\epsilon z)) \\
&=&\mathbb{P}(0\leq t_{0}+N_{0}(\eta +\epsilon z),\;M_{1}V\leq
t_{1}+N_{1}(\eta +\epsilon z)) \\
&=&\mathbb{P}(0\leq t_{0}^{\ast }+\epsilon N_{0}z,\;M_{1}(V+\mu )\leq
t_{1}^{\ast }+\epsilon N_{1}z) \\
&=&\mathbb{P}(0\leq \epsilon N_{0}z,\;M_{1}(V+\mu )\leq t_{1}^{\ast
}+\epsilon N_{1}z)
\end{eqnarray*}%
Since $t_{1}^{\ast }>0$, we can find a $t_{1}^{\ast \ast }$ such that $%
0<t_{1}^{\ast \ast }<t_{1}^{\ast }+\epsilon N_{1}z$ holds for every $%
\epsilon $ with $\left\vert \epsilon \right\vert $ small enough. If now $%
\epsilon >0$ then%
\begin{equation*}
f(x_{0}+\epsilon z)=\mathbb{P}(M_{1}(V+\mu )\leq t_{1}^{\ast }+\epsilon
N_{1}z)\geq \mathbb{P}(M_{1}(V+\mu )\leq t_{1}^{\ast \ast }).
\end{equation*}%
The r.h.s. in the above display is positive because $t_{1}^{\ast \ast }>0$
and because the distribution of $M_{1}(V+\mu )$ assigns positive mass to any
neighborhood of the origin, since the same is true for the distribution of $%
V+\mu $ and since $M_{1}$ maps neighborhoods of zero into neighborhoods of
zero. Setting $c=\mathbb{P}(M_{1}(V+\mu )\leq t_{1}^{\ast \ast })/2$, we
have $f(x_{0}+\epsilon z)\geq c>0$ for each sufficiently small $\epsilon >0$%
. Furthermore, for $\epsilon <0$ we have $f(x_{0}+\epsilon z)=0$, since $%
f(x_{0}+\epsilon z)\leq \mathbb{P}(0\leq \epsilon N_{0}z)=0$ in view of $%
N_{0}z>0$.{\hfill $\Box $}

\begin{lemma}
\label{linsert} Let $Z$ be a random vector with values in $\mathbf{R}^{p}$, $%
p\geq 1$, with a distribution that is absolutely continuous with respect to
Lebesgue measure on $\mathbf{R}^{p}$. Let $B$ be a $k\times p$ matrix, $%
k\geq 1$. Then the cdf $\mathbb{P(}BZ\leq \cdot )$\ of $BZ$, is
discontinuous at $t\in \mathbf{R}^{k}$ if and only if $\mathbb{P(}BZ\leq
t)>0 $ and if for some $i_{0}$, $1\leq i_{0}\leq k$, the $i_{0}$-th row of $%
B $ and the $i_{0}$-th component of $t$ are both zero, i.e., $B_{i_{0},\cdot
}=(0,\dots ,0)$ and $t_{i_{0}}=0$.
\end{lemma}

\textbf{Proof of Lemma \ref{linsert}:} To establish sufficiency of the above
condition, let $\mathbb{P(}BZ\leq t)>0$, $t_{i_{0}}=0$ and $B_{i_{0},\cdot
}=(0,\dots ,0)$ for some $i_{0}$, $1\leq i_{0}\leq k$. Then, of course, $%
\mathbb{P(}B_{i_{0},\cdot }Z=0)=1$. For $t_{n}=t-n^{-1}e_{i_{0}}$, where $%
e_{i_{0}}$ denotes the $i_{0}$-th unit vector in $\mathbf{R}^{k}$, we have $%
\mathbb{P(}BZ\leq t_{n})\leq \mathbb{P(}B_{i_{0},\cdot }Z\leq t_{n,i_{0}})=%
\mathbb{P(}B_{i_{0},\cdot }Z\leq -1/n)=0$ for every $n$. Consequently, $%
\mathbb{P(}BZ\leq t)$ is discontinuous at $t$.

To establish necessity, we first show the following: If $t_{n}\in \mathbf{R}%
^{k}$ is a sequence converging to $t\in \mathbf{R}^{k}$ as $n\rightarrow
\infty $, then every accumulation point of the sequence $\mathbb{P(}BZ\leq
t_{n})$ has the form 
\begin{equation}
\mathbb{P(}B_{i_{1},\cdot }Z\leq t_{i_{1}},\dots ,B_{i_{m},\cdot }Z\leq
t_{i_{m}},B_{i_{m+1},\cdot }Z<t_{i_{m+1}},\dots ,B_{i_{k},\cdot }Z<t_{i_{k}})
\label{insert1}
\end{equation}%
for some $m$, $0\leq m\leq k$, and for some permutation $(i_{1},\dots
,i_{k}) $ of $(1,\dots ,k)$. This can be seen as follows: Let $\alpha $ be
an accumulation point of $\mathbb{P(}BZ\leq t_{n})$. Then we may find a
subsequence such that $\mathbb{P(}BZ\leq t_{n})$ converges to $\alpha $
along this subsequence. From this subsequence we may even extract a further
subsequence along which each component of the $k\times 1$ vector $t_{n}$
converges to the corresponding component of $t$ monotonously, that is,
either from above or from below. Without loss of generality, we may also
assume that those components which converge from below are strictly
increasing. The resulting subsequence will be denoted by $n_{j}$ in the
sequel. Assume that the components of $t_{n_{j}}$ with indices $i_{1},\dots
,i_{m}$ converge from above, while the components with indices $%
i_{m+1},\dots ,i_{k}$ converge from below. Now 
\begin{equation}
\mathbb{P(}BZ\leq t_{n_{j}})=\int_{z\in \mathbf{R}^{k}}\prod_{s=1}^{k}%
\mathbf{1}_{(-\infty ,t_{n_{j},s}]}(z_{s})\mathbb{P}_{BZ}(dz),
\label{insert2}
\end{equation}%
where $\mathbb{P}_{BZ}$ denotes the distribution of $BZ$. The integrand in (%
\ref{insert2}) now converges to $\prod_{l=1}^{m}1_{(-\infty
,t_{i_{l}}]}(z_{i_{l}})\prod_{l=m+1}^{k}1_{(-\infty ,t_{i_{l}})}(z_{i_{l}})$
for all $z\in \mathbf{R}^{k}$ as $n_{j}\rightarrow \infty $. The r.h.s. of (%
\ref{insert2}) converges to the expression in (\ref{insert1}) as $%
n_{j}\rightarrow \infty $ by the Dominated Convergence Theorem, while the
l.h.s. of (\ref{insert2}) converges to $\alpha $ by construction. This
establishes the claim regarding (\ref{insert1}).

Now suppose that $\mathbb{P(}BZ\leq t)$ is discontinuous at $t$; i.e., there
exists a sequence $t_{n}$ converging to $t$ as $n\rightarrow \infty $, such
that $\mathbb{P(}BZ\leq t_{n})$ does not converge to $\mathbb{P(}BZ\leq t)$
as $n\rightarrow \infty $. From the sequence $t_{n}$ we can extract a
subsequence $t_{n_{s}}$ along which $\mathbb{P(}BZ\leq t_{n_{s}})$ converges
to a limit different from $\mathbb{P(}BZ\leq t)$ as $n_{s}\rightarrow \infty 
$. As shown above, the limit has to be of the form (\ref{insert1}) and $m<k$
has to hold. Consequently, the limit of $\mathbb{P(}BZ\leq t_{n_{s}})$ is
smaller than $\mathbb{P(}BZ\leq t)=\mathbb{P(}B_{i,\cdot }Z\leq
t_{i},\,i=1,\dots ,k)$. The difference of $\mathbb{P(}BZ\leq t)$ and the
limit of $\mathbb{P(}BZ\leq t_{n_{s}})$ is positive and because of (\ref%
{insert1}) can be written as%
\begin{equation*}
\mathbb{P(}B_{i_{j},\cdot }Z\leq t_{i_{j}}\text{ for each }j=1,\dots ,k\text{%
, }B_{i_{j},\cdot }Z=t_{i_{j}}\text{for some }j=m+1,\dots ,k)>0.
\end{equation*}

We thus see that $\mathbb{P(}B_{i_{j_{0}},\cdot }Z=t_{i_{j_{0}}})>0$ for
some ${j_{0}}$ satisfying $m+1\leq j_{0}\leq k$. As $Z$ is absolutely
continuous with respect to Lebesgue measure on $\mathbf{R}^{p}$, this can
only happen if $B_{i_{j_{0}},\cdot }=(0,\dots ,0)$ and $t_{i_{j_{0}}}=0$.{%
\hfill $\Box $}

\begin{lemma}
\label{lb.1} Suppose that $A\tilde{\theta}(q)$ and $\tilde{\theta}_{q}(q)$
are asymptotically correlated, i.e., $C_{\infty }^{(q)}\neq 0$, for some $q$
satisfying $\O <{q}\leq P$, and let $q^{\ast }$ denote the largest $q$ with
this property. Moreover let $\theta \in M_{q{^{\ast }}-1}$, let $\sigma $
satisfy $0<\sigma <\infty $, and let $t\in \mathbf{R}^{k}$. Then $G_{\infty
,\theta ,\sigma ,\gamma }(t)$\ is non-constant as a function of $\gamma \in
M_{q{^{\ast }}}\backslash M_{q{^{\ast }-1}}$. More precisely, there exist $%
\delta _{0}>0$ and $\rho _{0}$, $0<\rho _{0}<\infty $, such that 
\begin{equation}
\sup_{\overset{\gamma ^{(1)},\gamma ^{(2)}\in M_{q^{\ast }}\backslash M_{q{%
^{\ast }-1}}}{||\gamma ^{(i)}||<\rho _{0},i=1,2}}\left\vert G_{\infty
,\theta ,\sigma ,\gamma ^{(1)}}(t)-G_{\infty ,\theta ,\sigma ,\gamma
^{(2)}}(t)\right\vert >2\delta _{0}  \label{lb.1-1}
\end{equation}%
holds. The constants $\delta _{0}$ and $\rho _{0}$ can be chosen in such a
way that they depend only on $t$, $Q$, $A$, $\sigma $, and the critical
values $c_{p}$ for $\O <{p}\leq P$.
\end{lemma}

\begin{lemma}
\label{lb.2} Suppose that $A\tilde{\theta}(q)$ and $\tilde{\theta}_{q}(q)$
are asymptotically correlated, i.e., $C_{\infty }^{(q)}\neq 0$, for some $q$
satisfying $\O <{q}\leq P$, and let $q^{\ast }$ denote the largest $q$ with
this property. Suppose further that for some $p_{\odot }$ satisfying $\O %
\leq p_{\odot }<q^{\ast }$ either $p_{\odot }=0$ holds or that$\ p_{\odot
}>0 $ and $A[p_{\odot }]$ has a row of zeros. Then, for every $\theta \in
M_{p_{\odot }}$, every $\sigma $, $0<\sigma <\infty $, and every $t\in 
\mathbf{R}^{k}$ the quantity $G_{\infty ,\theta ,\sigma ,\gamma }(t)$ is
discontinuous as a function of $\gamma \in M_{q^{\ast }}$. More precisely,
for each $s=\O ,\dots ,p_{\odot }$, there exist vectors $\beta _{\ast }$ and 
$\gamma _{\ast }$ in $M_{q^{\ast }}$ and constants $\delta _{\ast }>0$ and $%
\epsilon _{\ast }>0$ such that 
\begin{equation}
\left\vert G_{\infty ,\theta ,\sigma ,\beta _{\ast }+\epsilon \gamma _{\ast
}}(t)-G_{\infty ,\theta ,\sigma ,\beta _{\ast }-\epsilon \gamma _{\ast
}}(t)\right\vert \geq \delta _{\ast }  \label{lb.2-1}
\end{equation}%
holds for every $\theta $ satisfying $\max \{p_{0}(\theta ),\O \}=s$ and for
every $\epsilon $ with $0<\epsilon <\epsilon _{\ast }$. The quantities $%
\delta _{\ast }$, $\epsilon _{\ast }$, $\beta _{\ast }$, and $\gamma _{\ast
} $ can be chosen in such a way that -- besides $t$, $Q$, $A$, $\sigma $,
and the critical values $c_{p}$ for $\O <{p}\leq P$ -- they depend on $%
\theta $ only through $\max \{p_{0}(\theta ),\O \}$.
\end{lemma}

Before we prove the above lemmas, we provide a representation of $G_{\infty
,\theta ,\sigma ,\gamma }(t)$ that will be useful in the following: For $%
0<p\leq P$ define $Z_{p}=\sum_{r=1}^{p}\xi _{\infty ,r}^{-2}C_{\infty
}^{(r)}W_{r}$, where $C_{\infty }^{(r)}$ has been defined after (\ref%
{zetainfinity}) and the random variables $W_{r}$ are independent normally
distributed with mean zero and variances $\sigma ^{2}\xi _{\infty ,r}^{2}$;
for convenience, let $Z_{0}$ denote the zero vector in $\mathbf{R}^{k}$.
Observe that $Z_{p}$, $p>0$, is normally distributed with mean zero and
variance-covariance matrix $\sigma ^{2}A[p]Q[p:p]^{-1}A[p]^{\prime }$ since
it has been shown in the proof of Proposition~4.4 in Leeb and P\"{o}tscher
(2006b) that the asymptotic variance-covariance matrix $\sigma
^{2}A[p]Q[p:p]^{-1}A[p]^{\prime }$ of $\sqrt{n}A\tilde{\theta}(p)$ can be
expressed as $\sum_{r=1}^{p}\sigma ^{2}\xi _{\infty ,r}^{-2}C_{\infty
}^{(r)}C_{\infty }^{(r)\prime }$. Also the joint distribution of $Z_{p}$ and
the set of variables $W_{r}$, $1\leq r\leq P$, is normal, with the
covariance vector between $Z_{p}$ and $W_{r}$ given by $\sigma ^{2}C_{\infty
}^{(r)}$ in case $r\leq p$; otherwise $Z_{p}$ and $W_{r}$ are independent.
Define the constants $\nu _{r}=\gamma _{r}+(Q[r:r]^{-1}Q[r:\lnot r]\gamma
\lbrack \lnot r])_{r}$ for $0<r\leq P$. It is now easy to see that for $%
p\geq p_{\ast }=\max \{p_{0}(\theta ),{\O }\}$ the quantity $\beta ^{(p)}$
defined in Proposition~\ref{p3.1} equals $-\sum_{r=p+1}^{P}\xi _{\infty
,r}^{-2}C_{\infty }^{(r)}\nu _{r}$. [This is seen as follows: It was noted
in Proposition~\ref{p3.1} that $\beta ^{(p)}=\lim_{n\rightarrow \infty }%
\sqrt{n}A(\eta _{n}(p)-\theta -\gamma /\sqrt{n})$ for $p\geq p_{0}(\theta )$%
, when $\eta _{n}(p)$ is defined as in (\ref{3.2}), but with $\theta +\gamma
/\sqrt{n}$ replacing $\theta $. Using the representation (20) of Leeb (2005)
and taking limits, the result follows if we observe that $\sqrt{n}\eta
_{n,r}(r)\longrightarrow \nu _{r}$ for $r>p\geq p_{0}(\theta )$.] The cdf in
(\ref{3.7}) can now be written as%
\begin{align}
\mathbb{P}\left( Z_{p_{\ast }}\leq t+\sum_{r=p_{\ast }+1}^{P}\xi _{\infty
,r}^{-2}C_{\infty }^{(r)}\nu _{r}\right) \prod_{q=p_{\ast }+1}^{P}\mathbb{P}%
(|W_{q}+\nu _{q}|& <c_{q}\sigma \xi _{\infty ,q})  \notag \\
+\sum_{p=p_{\ast }+1}^{P}\mathbb{P}\left( Z_{p}\leq t+\sum_{r=p+1}^{P}\xi
_{\infty ,r}^{-2}C_{\infty }^{(r)}\nu _{r},\;|W_{p}+\nu _{p}|\geq
c_{p}\sigma \xi _{\infty ,p}\right) \prod_{q=p+1}^{P}\mathbb{P}(|W_{q}+\nu
_{q}|& <c_{q}\sigma \xi _{\infty ,q}).  \label{ab.1}
\end{align}%
That the terms corresponding to $p=p_{\ast }$ in (\ref{ab.1}) and (\ref{3.7}%
) agree is obvious. Furthermore, for each $p>p_{\ast }$ the terms under the
product sign in (\ref{ab.1}) and (\ref{3.7}) coincide by definition of the
function $\Delta _{s}(a,b)$. It is also easy to see that the conditional
distribution of $W_{p}$ given $Z_{p}=z$ is Gaussian with mean $b_{\infty
,p}z $ and variance $\sigma ^{2}\zeta _{\infty ,p}^{2}$. Consequently, the
probability of the event $\{|W_{p}+\nu _{p}|\geq c_{p}\sigma \xi _{\infty
,p}\}$ conditional on $Z_{p}=z$ is given by the integrand shown in (\ref{3.7}%
). Since $Z_{p}$ has distribution $\Phi _{\infty ,p}$ as noted above, it
follows that (\ref{ab.1}) and (\ref{3.7}) agree.

\begin{remark}
\normalfont\label{rb.3}If $C_{\infty }^{(p)}=0$ for $p>p_{\ast }$, then in
view of the above discussion $Z_{p_{\ast }}=Z_{p}=Z_{P}$, and hence $\Phi
_{\infty ,p_{\ast }}=\Phi _{\infty ,p}=\Phi _{\infty ,P}$, holds for all $%
p>p_{\ast }$. Using the independence of $W_{r}$, $r>p_{\ast }$, from $%
Z_{p_{\ast }}$, inspection of (\ref{ab.1}) shows that $G_{\infty ,\theta
,\sigma ,\gamma }$ reduces to $\Phi _{\infty ,P}$;\ see also Leeb (2006,
Remark 5.2).
\end{remark}

\textbf{Proof of Lemma \ref{lb.1}:} From (\ref{ab.1}) (or (\ref{3.7})) it
follows that the map $\gamma \mapsto G_{\infty ,\theta ,\sigma ,\gamma }(t)$
depends only on $t$, $Q$, $A$, $\sigma $, the critical values $c_{p}$ for $%
\O <{p}\leq P$, as well as on $\theta $; however, the dependence on $\theta $%
\ is only through $p_{\ast }=\max \{p_{0}(\theta ),\O \}$. It hence suffices
to find, for each possible value of $p_{\ast }$ in the range $p_{\ast }=\O %
,\dots ,q^{\ast }-1$, constants $0<\rho _{0}<\infty $ and $\delta _{0}>0$
such that (\ref{lb.1-1}) is satisfied for some (and hence all) $\theta $
returning this particular value of $p_{\ast }=\max \{p_{0}(\theta ),\O \}$.
For this in turn it is sufficient to show that for every $\theta \in
M_{q^{\ast }-1}$ the quantity $G_{\infty ,\theta ,\sigma ,\gamma }(t)$ is
non-constant as a function of $\gamma \in M_{q^{\ast }}\backslash M_{q{%
^{\ast }-1}}$.

Let $\theta \in M_{q^{\ast }-1}$ and assume that $G_{\infty ,\theta ,\sigma
,\gamma }(t)$ is constant in $\gamma \in M_{q^{\ast }}\backslash M_{q{^{\ast
}-1}}$. Observe that, by assumption, $C_{\infty }^{(q^{\ast })}$ is non-zero
while $C_{\infty }^{(p)}=0$ for $p>q^{\ast }$. For $\gamma \in M_{q^{\ast }}$%
, we clearly have $\nu _{q^{\ast }}=\gamma _{q^{\ast }}$ and $\nu _{r}=0$
for $r>q^{\ast }$. Letting $\gamma _{q^{\ast }-1}\rightarrow \infty $ while $%
\gamma _{q^{\ast }}$ is held fixed, we see that $\nu _{q^{\ast
}-1}\rightarrow \infty $; hence,%
\begin{equation*}
\mathbb{P}(|W_{q^{\ast }-1}+\nu _{q^{\ast }-1}|<c_{q^{\ast }-1}\sigma \xi
_{\infty ,q^{\ast }-1})\rightarrow 0.
\end{equation*}%
It follows that (\ref{ab.1}) converges to

\begin{gather}
\mathbb{P}\left( Z_{q^{\ast }-1}\leq t+\xi _{\infty ,q^{\ast
}}^{-2}C_{\infty }^{(q^{\ast })}\gamma _{q^{\ast }}\right) \mathbb{P}%
(|W_{q^{\ast }}+\gamma _{q^{\ast }}|<c_{q^{\ast }}\sigma \xi _{\infty
,q^{\ast }})\prod_{q=q^{\ast }+1}^{P}\mathbb{P}(|W_{q}|<c_{q}\sigma \xi
_{\infty ,q})  \notag \\
+\mathbb{P}\left( Z_{q^{\ast }}\leq t,\;|W_{q^{\ast }}+\gamma _{q^{\ast
}}|\geq c_{q^{\ast }}\sigma \xi _{\infty ,q^{\ast }}\right) \prod_{q=q^{\ast
}+1}^{P}\mathbb{P}(|W_{q}|<c_{q}\sigma \xi _{\infty ,q})  \label{ab.2} \\
+\sum_{p=q^{\ast }+1}^{P}\mathbb{P}\left( Z_{p}\leq t,\;|W_{p}|\geq
c_{p}\sigma \xi _{\infty ,p}\right) \prod_{q=p+1}^{P}\mathbb{P}%
(|W_{q}|<c_{q}\sigma \xi _{\infty ,q}).  \notag
\end{gather}%
By assumption, the expression in the above display is constant in $\gamma
_{q^{\ast }}\in \mathbf{R}\backslash \{0\}$. Dropping the terms that do not
depend on $\gamma _{q^{\ast }}$ and observing that $\mathbb{P}%
(|W_{q}|<c_{q}\sigma \xi _{\infty ,q})$ is never zero for $q>q^{\ast }>\O $,
we see that 
\begin{align}
\mathbb{P}\left( Z_{q^{\ast }-1}\leq t+\xi _{\infty ,q^{\ast
}}^{-2}C_{\infty }^{(q^{\ast })}\gamma _{q^{\ast }}\right) \mathbb{P}%
(|W_{q^{\ast }}+\gamma _{q^{\ast }}|& <c_{q^{\ast }}\sigma \xi _{\infty
,q^{\ast }})  \notag \\
& +\mathbb{P}\left( Z_{q^{\ast }}\leq t,\;|W_{q^{\ast }}+\gamma _{q^{\ast
}}|\geq c_{q^{\ast }}\sigma \xi _{\infty ,q^{\ast }}\right)  \label{ab.3}
\end{align}%
has to be constant in $\gamma _{q^{\ast }}\in \mathbf{R}\backslash \{0\}$.
We now show that the expression in (\ref{ab.3}) is in fact constant in $%
\gamma _{q^{\ast }}\in \mathbf{R}$: Observe first that $\mathbb{P}%
(|W_{q^{\ast }}+\gamma _{q^{\ast }}|<c_{q^{\ast }}\sigma \xi _{\infty
,q^{\ast }})$ is positive and continuous in $\gamma _{q^{\ast }}\in \mathbf{R%
}$; also the probability $\mathbb{P}\left( Z_{q^{\ast }}\leq t,\;|W_{q^{\ast
}}+\gamma _{q^{\ast }}|\geq c_{q^{\ast }}\sigma \xi _{\infty ,q^{\ast
}}\right) $ is continuous in $\gamma _{q^{\ast }}\in \mathbf{R}$ since $%
W_{q^{\ast }}$, being normal with mean zero and positive variance, is
absolutely continuously distributed. Concerning the remaining term in (\ref%
{ab.3}), we note that $Z_{q^{\ast }-1}=MV$ where $M=[\xi _{\infty
,1}^{-2}C_{\infty }^{(1)},\ldots ,\xi _{\infty ,q^{\ast }-1}^{-2}C_{\infty
}^{(q^{\ast }-1)}]$ and $V=(W_{1},\dots ,W_{q^{\ast }-1})^{\prime }$. In
case no row of $M$ is identically zero, Lemma \ref{linsert} shows that also $%
\mathbb{P}\left( Z_{q^{\ast }-1}\leq t+\xi _{\infty ,q^{\ast
}}^{-2}C_{\infty }^{(q^{\ast })}\gamma _{q^{\ast }}\right) $ is continuous
in $\gamma _{q^{\ast }}\in \mathbf{R}$. Hence, in this case (\ref{ab.3}) is
indeed constant for all $\gamma _{q^{\ast }}\in \mathbf{R}$. In case a row
of $M$ is identically zero, define $N=\xi _{\infty ,q^{\ast }}^{-2}C_{\infty
}^{(q^{\ast })}$ and rewrite the probability in question as $\mathbb{P}%
\left( MV\leq t+N\gamma _{q^{\ast }}\right) $. Note that $(M:N)$ has full
row-rank $k$, since%
\begin{equation}
(M:N)diag[\xi _{\infty ,1}^{2},\ldots ,\xi _{\infty ,q^{\ast
}}^{2}](M:N)^{\prime }=\sum_{r=1}^{q^{\ast }}\xi _{\infty ,r}^{-2}C_{\infty
}^{(r)}C_{\infty }^{(r)\prime }=\sum_{r=1}^{P}\xi _{\infty ,r}^{-2}C_{\infty
}^{(r)}C_{\infty }^{(r)\prime }=AQ^{-1}A^{\prime }  \label{matrix}
\end{equation}%
by definition of $q^{\ast }$ and since the latter matrix is non-singular in
view of $rank$ $A=k$. Lemma \ref{la.2} then shows that there exists a $%
\gamma _{q^{\ast }}^{(0)}\in \mathbf{R}$, $z\in \{-1,1\}$, and a constant $%
c>0$ such that $\mathbb{P}\left( MV\leq t+N(\gamma _{q^{\ast }}^{(0)}-\delta
z)\right) =0$ and $\mathbb{P}\left( MV\leq t+N(\gamma _{q^{\ast
}}^{(0)}+\delta z)\right) \geq c$ holds for arbitrary small $\delta >0$.
Observe that $\gamma _{q^{\ast }}^{(0)}-\delta z$ as well as $\gamma
_{q^{\ast }}^{(0)}-\delta z$ are non-zero for sufficiently small $\delta >0$%
. But then (\ref{ab.3}) -- being constant for $\gamma _{q^{\ast }}\in 
\mathbf{R}\backslash \{0\}$ -- gives the same value for $\gamma _{q^{\ast
}}=\gamma _{q^{\ast }}^{(0)}-\delta z$ and $\gamma _{q^{\ast }}=\gamma
_{q^{\ast }}^{(0)}+\delta z$ and all sufficiently small $\delta >0$. Letting 
$\delta $ go to zero in this equality and using the continuity properties
for the second and third probability in (\ref{ab.3}) noted above we obtain
that 
\begin{multline*}
c\mathbb{P}(|W_{q^{\ast }}+\gamma _{q^{\ast }}^{(0)}|<c_{q^{\ast }}\sigma
\xi _{\infty ,q^{\ast }})+\mathbb{P}\left( Z_{q^{\ast }}\leq t,\;|W_{q^{\ast
}}+\gamma _{q^{\ast }}^{(0)}|\geq c_{q^{\ast }}\sigma \xi _{\infty ,q^{\ast
}}\right) \\
\leq \liminf_{\delta \downarrow 0}\mathbb{P}\left( Z_{q^{\ast }-1}\leq t+\xi
_{\infty ,q^{\ast }}^{-2}C_{\infty }^{(q^{\ast })}(\gamma _{q^{\ast
}}^{(0)}+\delta z)\right) \mathbb{P}(|W_{q^{\ast }}+\gamma _{q^{\ast
}}^{(0)}|<c_{q^{\ast }}\sigma \xi _{\infty ,q^{\ast }}) \\
+\mathbb{P}\left( Z_{q^{\ast }}\leq t,\;|W_{q^{\ast }}+\gamma _{q^{\ast
}}^{(0)}|\geq c_{q^{\ast }}\sigma \xi _{\infty ,q^{\ast }}\right) \\
=\liminf_{\delta \downarrow 0}\mathbb{P}\left( Z_{q^{\ast }-1}\leq t+\xi
_{\infty ,q^{\ast }}^{-2}C_{\infty }^{(q^{\ast })}(\gamma _{q^{\ast
}}^{(0)}-\delta z)\right) \mathbb{P}(|W_{q^{\ast }}+\gamma _{q^{\ast
}}^{(0)}|<c_{q^{\ast }}\sigma \xi _{\infty ,q^{\ast }}) \\
+\mathbb{P}\left( Z_{q^{\ast }}\leq t,\;|W_{q^{\ast }}+\gamma _{q^{\ast
}}^{(0)}|\geq c_{q^{\ast }}\sigma \xi _{\infty ,q^{\ast }}\right) \\
=\mathbb{P}\left( Z_{q^{\ast }}\leq t,\;|W_{q^{\ast }}+\gamma _{q^{\ast
}}^{(0)}|\geq c_{q^{\ast }}\sigma \xi _{\infty ,q^{\ast }}\right)
\end{multline*}%
which is impossible since $c>0$ and $\mathbb{P}(|W_{q^{\ast }}+\gamma
_{q^{\ast }}^{(0)}|<c_{q^{\ast }}\sigma \xi _{\infty ,q^{\ast }})>0$. Hence
we have shown that (\ref{ab.3}) is indeed constant for all $\gamma _{q^{\ast
}}\in \mathbf{R}$.

Now write $Z$, $W$, $C$, $\tau $, and $x$ for $Z_{q^{\ast }-1}-t$, $%
-W_{q^{\ast }}/\sigma \xi _{\infty ,q^{\ast }}$, $\sigma \xi _{\infty
,q^{\ast }}^{-1}C_{\infty }^{(q^{\ast })}$, $c_{q^{\ast }}$, and $\gamma
_{q^{\ast }}/\sigma \xi _{\infty ,q^{\ast }}$, respectively. Upon observing
that $Z_{q^{\ast }}$ equals $Z_{q^{\ast }-1}+\xi _{\infty ,q^{\ast
}}^{-2}C_{\infty }^{(q^{\ast })}W_{q^{\ast }}$, it is easy to see that (\ref%
{ab.3}) can be written as in (\ref{la.1-1}). By our assumptions, this
expression is constant in $x=\gamma _{q^{\ast }}/\sigma \xi _{\infty
,q^{\ast }}\in \mathbf{R}$. Lemma~\ref{la.1} then entails that either $C=0$
or that $\mathbb{P}(Z\leq Cx)=0$ for each $x\in \mathbf{R}$. Since $C$
equals $\sigma \xi _{\infty ,q^{\ast }}^{-1}C_{\infty }^{(q^{\ast })}$, it
is non-zero by assumption. Hence,%
\begin{equation*}
\mathbb{P}\left( Z_{q^{\ast }-1}\leq t+\xi _{\infty ,q^{\ast
}}^{-2}C_{\infty }^{(q^{\ast })}\gamma _{q^{\ast }}\right) =0
\end{equation*}%
must hold for every value of $\gamma _{q^{\ast }}$. But the above
probability is just the conditional probability that $Z_{q^{\ast }}\leq t$
given $W_{q^{\ast }}=-\gamma _{q^{\ast }}$. It follows that $\mathbb{P}%
(Z_{q^{\ast }}\leq t)$ equals zero as well. By our assumption $C_{\infty
}^{(p)}=0$ for $p>q^{\ast }$, and hence $Z_{q^{\ast }}=Z_{P}$. We thus
obtain $\mathbb{P}(Z_{P}\leq t)=0$, a contradiction with the fact that $%
Z_{P} $ is a Gaussian random variable on $\mathbf{R}^{k}$ with non-singular
variance-covariance matrix $\sigma ^{2}AQ^{-1}A^{\prime }$. {\hfill $\Box $}

Inspection of the above proof shows that it can be simplified if the claim
of non-constancy of $G_{\infty ,\theta ,\sigma ,\gamma }(t)$ as a function
of $\gamma \in M_{q{^{\ast }}}\backslash M_{q{^{\ast }-1}}$ in Lemma \ref%
{lb.1}\textbf{\ }is weakened to non-constancy for $\gamma \in M_{q{^{\ast }}%
} $. The strong form of the lemma as given here is needed in the proof of
Proposition \ref{p.ae1}.

\textbf{Proof of Lemma \ref{lb.2}:} Let $p_{\oplus }$ be the largest index $%
p $, $\O \leq p\leq P$, for which $A[p]$ has a row of zeroes, and set $%
p_{\oplus }=0$ if no such index exists. We first show that $p_{\oplus }$
satisfies $p_{\oplus }<q^{\ast }$. Suppose $p_{\oplus }\geq q^{\ast }$ would
hold. Since $Z_{p_{\oplus }}$ is a Gaussian random vector with mean zero and
variance-covariance matrix $\sigma ^{2}A[p_{\oplus }]Q[p_{\oplus }:p_{\oplus
}]^{-1}A[p_{\oplus }]^{\prime }$, at least one component of $Z_{p_{\oplus }}$
is equal to zero with probability one. However, $Z_{p_{\oplus }}$ equals $%
Z_{P}$ because of $p_{\oplus }\geq q^{\ast }$ and the definition of $q^{\ast
}$. This leads to a contradiction since $Z_{P}$ has the non-singular
variance-covariance matrix $\sigma ^{2}AQ^{-1}A^{\prime }$. Without loss of
generality, we may hence assume that $p_{\odot }=p_{\oplus }$.

In view of the discussion in the first paragraph of the proof of Lemma \ref%
{lb.1}, it suffices to establish, for each possible value $s$ in the range $%
\O \leq s\leq p_{\odot }$, the result (\ref{lb.2-1}) for some $\theta $ with 
$s=\max \{p_{0}(\theta ),\O \}=p_{\ast }$. Now fix such an $s$ and $\theta $
(as well as, of course, $t$, $Q$, $A$, $\sigma $, and the critical values $%
c_{p}$ for $\O <{p}\leq P$). Then (\ref{ab.1}) expresses the map $\gamma
\mapsto G_{\infty ,\theta ,\sigma ,\gamma }(t)$ in terms of $\nu =(\nu
_{1},\dots ,\nu _{P})^{\prime }$. It is easy to see that the correspondence
between $\gamma $ and $\nu $ is a linear bijection from $\mathbb{R}^{P}$
onto itself, and that $\gamma \in M_{q^{\ast }}$ if and only if $\nu \in
M_{q^{\ast }}$. It is hence sufficient to find a $\delta _{\ast }>0$ and
vectors $\nu $ and $\mu $ in $M_{q^{\ast }}$ such that (\ref{ab.1}) with $%
\nu +\epsilon \mu $ in place of $\nu $ and (\ref{ab.1}) with $\nu -\epsilon
\mu $ in place of $\nu $ differ by at least $\delta _{\ast }$ for
sufficiently small $\epsilon >0$. Note that (\ref{ab.1}) is the sum of $%
P-p_{\ast }+1$ terms indexed by $p=p_{\ast },\dots ,P$. We shall now show
that $\nu $ and $\mu $ can be chosen in such a way that, when replacing $\nu 
$ with $\nu +\epsilon \mu $ and $\nu -\epsilon \mu $, respectively, (i) the
resulting terms in (\ref{ab.1}) corresponding to $p=p_{\odot }$ differ by
some $d>0$, while (ii) the difference of the other terms becomes arbitrarily
small, provided that $\epsilon >0$ is sufficiently small.

Consider first the case where $s=p_{\ast }=p_{\odot }$. Using the shorthand
notation%
\begin{equation*}
g(\nu )=\mathbb{P}\left( Z_{p_{\odot }}\leq t+\sum_{r=p_{\odot }+1}^{q^{\ast
}}\xi _{\infty ,r}^{-2}C_{\infty }^{(r)}\nu _{r}\right) ,
\end{equation*}%
note that the $p_{\odot }$-th term in (\ref{ab.1}) is given by $g(\nu )$
multiplied by a product of positive probabilities which are continuous in $%
\nu $. To prove property (i) it thus suffices to find a constant $c>0$, and
vectors $\nu $ and $\mu $ in $M_{q^{\ast }}$ such that $|g(\nu +\epsilon \mu
)-g(\nu -\epsilon \mu )|\geq c$ holds for each sufficiently small $\epsilon
>0$.

In the sub-case $p_{\odot }=0$ choose $c=1$, set 
\begin{equation*}
\nu =-[C_{\infty }^{(1)},\ldots ,C_{\infty }^{(P)}]^{\prime }\left[
\sum_{r=1}^{P}\xi _{\infty ,r}^{-2}C_{\infty }^{(r)}C_{\infty }^{(r)\prime }%
\right] ^{-1}t
\end{equation*}%
and%
\begin{equation*}
\mu =[C_{\infty }^{(1)},\ldots ,C_{\infty }^{(P)}]^{\prime }\left[
\sum_{r=1}^{P}\xi _{\infty ,r}^{-2}C_{\infty }^{(r)}C_{\infty }^{(r)\prime }%
\right] ^{-1}(1,\ldots ,1)^{\prime },
\end{equation*}%
observing that the matrix to be inverted is indeed non-singular, since -- as
discussed after Lemma \ref{lb.2} -- it is up to a multiplicative factor $%
\sigma ^{2}$\ identical to the variance-covariance matrix $\sigma
^{2}AQ^{-1}A^{\prime }$\ of $Z_{P}$. But then $\nu $ and $\mu $ satisfy $%
\sum_{r=p_{\odot }+1}^{q^{\ast }}\xi _{\infty ,r}^{-2}C_{\infty }^{(r)}\nu
_{r}=-t$ and $\sum_{r=p_{\odot }+1}^{q^{\ast }}\xi _{\infty
,r}^{-2}C_{\infty }^{(r)}\mu _{r}=(1,\dots ,1)^{\prime }$ if we note that by
the definition of $q^{\ast }$%
\begin{equation*}
\sum_{r=p_{\odot }+1}^{q^{\ast }}\xi _{\infty ,r}^{-2}C_{\infty }^{(r)}\nu
_{r}=\sum_{r=1}^{P}\xi _{\infty ,r}^{-2}C_{\infty }^{(r)}\nu _{r}
\end{equation*}%
holds and that a similar relation holds with $\mu $ replacing $\nu $. Since $%
Z_{p_{\odot }}=Z_{0}=0\in \mathbb{R}^{k}$, it is then obvious that $g(\nu
+\epsilon \mu )$ and $g(\nu -\epsilon \mu )$ differ by $1$ for each $%
\epsilon >0$.

In the other sub-case $p_{\odot }>0$, define $M=[\xi _{\infty
,1}^{-2}C_{\infty }^{(1)},\ldots ,\xi _{\infty ,p_{\odot }}^{-2}C_{\infty
}^{(p_{\odot })}]$, $N=[\xi _{\infty ,p_{\odot }+1}^{-2}C_{\infty
}^{(p_{\odot }+1)},\ldots ,\xi _{\infty ,q^{\ast }}^{-2}C_{\infty
}^{(q^{\ast })}]$, and $V=(W_{1},\dots ,W_{p_{\odot }})^{\prime }$. It is
then easy to see that $g(\nu )$ equals $f((\nu _{p_{\odot }+1},\dots ,\nu
_{q^{\ast }})^{\prime })$, with $f$ defined as in Lemma~\ref{la.2}, and that 
$M$ has a row of zeros. Furthermore, the matrix $(M:N)$ has rank $k$ by the
same argument as in the proof of Lemma \ref{lb.1}; cf. (\ref{matrix}). By
Lemma~\ref{la.2}, we thus obtain vectors $x_{0}$ and $z$, and a $c>0$ such
that $|f(x_{0}+\epsilon z)-f(x_{0}-\epsilon z)|\geq c$ holds for each
sufficiently small $\epsilon >0$. Setting $(\nu _{p_{\odot }+1},\dots ,\nu
_{q^{\ast }})^{\prime }=x_{0}$, $(\mu _{p_{\odot }+1},\dots ,\mu _{q^{\ast
}})^{\prime }=z$, setting $\nu \lbrack \lnot q^{\ast }]$, and $\mu \lbrack
\lnot q^{\ast }]$ each equal to zero, and setting $\nu \lbrack p_{\odot }]$
and $\mu \lbrack p_{\odot }]$ to arbitrary values, we see that $g(\nu \pm
\epsilon \mu )$ has the desired properties.

To complete the proof in case $s=p_{\ast }=p_{\odot }$, we need to establish
property (ii) for which it suffices to show that, for $p>p_{\odot }$, the $p$%
-th term in (\ref{ab.1}) depends continuously on $\nu $. For $p>q^{\ast }$,
the $p$-th term does not depend on $\nu $, because $C_{\infty }^{(r)}=0$ for 
$r=q^{\ast },\dots ,P$. For $p$ satisfying $p_{\odot }<p\leq q^{\ast }$, it
suffices to show that%
\begin{equation*}
h(\nu _{p},\ldots ,\nu _{q^{\ast }})=\mathbb{P}\left( Z_{p}\leq
t+\sum_{r=p+1}^{q^{\ast }}\xi _{\infty ,r}^{-2}C_{\infty }^{(r)}\nu
_{r},\;|W_{p}+\nu _{p}|\geq c_{p}\sigma \xi _{\infty ,p}\right)
\end{equation*}%
is a continuous function. Suppose that $(\nu _{p}^{(m)},\ldots ,\nu
_{q^{\ast }}^{(m)})$ converges to $(\nu _{p},\ldots ,\nu _{q^{\ast }})$ as $%
m\rightarrow \infty $. For arbitrary $\alpha >0$, $\sum_{r=p+1}^{q^{\ast
}}\xi _{\infty ,r}^{-2}C_{\infty }^{(r)}\nu _{r}$ and $\sum_{r=p+1}^{q^{\ast
}}\xi _{\infty ,r}^{-2}C_{\infty }^{(r)}\nu _{r}^{(m)}$ differ by less than $%
\alpha $ in each coordinate, provided that $m$ is sufficiently large. This
implies%
\begin{align*}
& \limsup_{m\rightarrow \infty }h(\nu _{p}^{(m)},\ldots ,\nu _{q^{\ast
}}^{(m)}) \\
& \leq \limsup_{m\rightarrow \infty }\mathbb{P(}Z_{p}\leq
t+\sum_{r=p+1}^{q^{\ast }}\xi _{\infty ,r}^{-2}C_{\infty }^{(r)}\nu
_{r}+\alpha (1,\ldots ,1)^{\prime },\;|W_{p}+\nu _{p}^{(m)}|\geq c_{p}\sigma
\xi _{\infty ,p}) \\
& =\mathbb{P(}Z_{p}\leq t+\sum_{r=p+1}^{q^{\ast }}\xi _{\infty
,r}^{-2}C_{\infty }^{(r)}\nu _{r}+\alpha (1,\ldots ,1)^{\prime
},\;|W_{p}+\nu _{p}|\geq c_{p}\sigma \xi _{\infty ,p}),
\end{align*}%
observing that the latter probability is obviously continuous in the single
variable $\nu _{p}$ (since $W_{p}$ has an absolutely continuous
distribution). Letting $\alpha $ decrease to zero we obtain $%
\limsup_{m\rightarrow \infty }h(\nu _{p}^{(m)},\ldots ,\nu _{q^{\ast
}}^{(m)})\leq h(\nu _{p},\ldots ,\nu _{q^{\ast }})$. A similar argument
establishes $\liminf_{m\rightarrow \infty }h(\nu _{p}^{(m)},\ldots ,\nu
_{q^{\ast }}^{(m)})\geq \mathbb{P(}Z_{p}<t+\sum_{r=p+1}^{q^{\ast }}\xi
_{\infty ,r}^{-2}C_{\infty }^{(r)}\nu _{r},\;|W_{p}+\nu _{p}|\geq
c_{p}\sigma \xi _{\infty ,p})$. The proof of the continuity of $h$ is then
complete if we can show that $\mathbb{P}\left( Z_{p}\leq \cdot ,\;|W_{p}+\nu
_{p}|\geq c_{p}\sigma \xi _{\infty ,p}\right) $ is continuous or,
equivalently, that $\mathbb{P}\left( Z_{p}\leq \cdot \left\vert |W_{p}+\nu
_{p}|\geq c_{p}\sigma \xi _{\infty ,p}\right. \right) $ is a continuous cdf.
Since $p>p_{\odot }$, the variance-covariance matrix $\sigma
^{2}A[p]Q[p:p]^{-1}A[p]^{\prime }$ of $Z_{p}$\ does only have non-zero
diagonal elements. Consequently, when representing $Z_{p}$ as $%
B(W_{1},\ldots ,W_{p})^{\prime }$, the matrix $B$ cannot have rows that
consist entirely of zeros. The conditional distribution of $(W_{1},\ldots
,W_{p})^{\prime }$ given the event $\{|W_{p}+\nu _{p}|\geq c_{p}\sigma \xi
_{\infty ,p}\}$ is clearly absolutely continuous w.r.t. $p$-dimensional
Lebesgue measure. But then Lemma \ref{linsert} delivers the desired result.

The case where $s=p_{\ast }<p_{\odot }$ is reduced to the previously
discussed case as follows: It is easy to see that, for $\nu _{p_{\odot
}}\rightarrow \infty $, the expression in (\ref{ab.1}) converges to a limit
uniformly w.r.t. all $\nu _{p}$ with $p\neq p_{\odot }$. Then observe that
this limit is again of the form (\ref{ab.1}) but now with $p_{\odot }$
taking the r\^{o}le of $p_{\ast }$. {\hfill $\Box $}

\section{Non-Uniformity of the Convergence of the Finite-Sample Cdf to the
Large-Sample Limit\label{ae}}

\begin{proposition}
\label{p.ae1}

\begin{enumerate}
\item \label{p.ae1.a} Suppose that $A\tilde{\theta}(q)$ and $\tilde{\theta}%
_{q}(q)$ are asymptotically correlated, i.e., $C_{\infty }^{(q)}\neq 0$, for
some $q$ satisfying $\O <{q}\leq P$, and let $q^{\ast }$ denote the largest $%
q$ with this property. Then for every $\theta \in M_{q{^{\ast }}-1}$, every $%
\sigma $, $0<\sigma <\infty $, and every $t\in \mathbf{R}^{k}$ there exists
a $\rho $, $0<\rho <\infty $, such that%
\begin{equation}
\liminf_{n\rightarrow \infty }\sup_{\overset{\vartheta \in M_{q^{\ast }}}{%
||\vartheta -\theta ||<\rho /\sqrt{n}}}\left\vert G_{n,\vartheta ,\sigma
}(t)-G_{\infty ,\vartheta ,\sigma }(t)\right\vert >0  \label{p3.2-1}
\end{equation}%
holds. The constant $\rho $ may be chosen in such a way that it depends only
on $t$, $Q$, $A$, $\sigma $, and the critical values $c_{p}$ for $\O <{p}%
\leq P$.

\item \label{p.ae1.b} Suppose that $A\tilde{\theta}(q)$ and $\tilde{\theta}%
_{q}(q)$ are asymptotically uncorrelated, i.e., $C_{\infty }^{(q)}=0$, for
all $q$ satisfying $\O <{q}\leq P$. Then $G_{n,\theta ,\sigma }$ converges
to $\Phi _{\infty ,P}$ in total variation uniformly in $\theta \in \mathbf{R}%
^{P}$; more precisely%
\begin{equation*}
\sup_{\theta \in \mathbf{R}^{P}}\sup_{\overset{\sigma \in \mathbf{R}}{\sigma
_{\ast }\leq \sigma \leq \sigma ^{\ast }}}\left\vert \left\vert G_{n,\theta
,\sigma }-\Phi _{\infty ,P}\right\vert \right\vert _{TV}\overset{%
n\rightarrow \infty }{\longrightarrow }0
\end{equation*}%
holds for any constants $\sigma _{\ast }$ and $\sigma ^{\ast }$ satisfying $%
0<\sigma _{\ast }\leq \sigma ^{\ast }<\infty $.
\end{enumerate}
\end{proposition}

Under the assumptions of Proposition~\ref{p.ae1}(\ref{p.ae1.a}), we see that
convergence of $G_{n,\theta ,\sigma }(t)$ to $G_{\infty ,\theta ,\sigma }(t)$
is non-uniform over shrinking `tubes' around $M_{q^{\ast }-1}$ that are
contained in $M_{q^{\ast }}$. [On the complement of a tube with a \emph{fixed%
} positive radius, i.e., on the set $U=\{\theta \in \mathbf{R}^{P}:\,|\theta
_{q^{\ast }}|\geq r\}$ with fixed $r>0$, convergence of $G_{n,\theta ,\sigma
}(t)$ to $G_{\infty ,\theta ,\sigma }(t)$ is in fact uniform (even with
respect to the total variation distance), as can be shown. Note that for $%
\theta \in U$ the cdf $G_{\infty ,\theta ,\sigma }(t)$ reduces to the
Gaussian cdf $\Phi _{\infty ,P}(t)$, i.e., to the asymptotic distribution of
the least-squares estimator based on the overall model; cf. Remark~\ref{rb.3}%
.] A precursor to Proposition~\ref{p.ae1}(\ref{p.ae1.a}) is Corollary~5.5 of
Leeb and P\"{o}tscher (2003) which establishes (\ref{p3.2-1}) in the special
case where $\O =0$ and where $A$ is the $P\times P$ identity matrix.
Proposition \ref{p.ae1}(\ref{p.ae1.b}) describes an exceptional case where
convergence is uniform. [In this case $G_{\infty ,\theta ,\sigma }$ reduces
to the Gaussian cdf $\Phi _{\infty ,P}$ for all $\theta $ and $\Phi _{\infty
,P}=\Phi _{\infty ,p}$, $\O \leq p\leq P$, holds; cf. Remark \ref{rb.3}.]
Recall that under the assumptions of part (\ref{p.ae1.b}) of Proposition \ref%
{p.ae1} we necessarily always have (i) $\O >0$, and (ii) $rank$ $A[\O ]=k$;
cf. Proposition~4.4 in Leeb and P\"{o}tscher (2006b).

\textbf{Proof of Proposition \ref{p.ae1}:} We first prove part (\ref{p.ae1.a}%
). As noted at the beginning of the proof of Lemma \ref{lb.1}, the map $%
\gamma \mapsto G_{\infty ,\theta ,\sigma ,\gamma }(t)$ depends only on $t$, $%
Q$, $A$, $\sigma $, the critical values $c_{p}$ for $\O <{p}\leq P$, as well
as on $\theta $, but the dependence on $\theta $\ is only through $p_{\ast
}=\max \{p_{0}(\theta ),\O \}$. It hence suffices to find, for each possible
value of $p_{\ast }$ in the range $p_{\ast }=\O ,\dots ,q^{\ast }-1$, a
constant $0<\rho <\infty $ such that (\ref{p3.2-1}) is satisfied for some
(and hence all) $\theta $ returning this particular value of $p_{\ast }=\max
\{p_{0}(\theta ),\O \}$. For this in turn it is sufficient to show that
given such a $\theta $ we can find a $\gamma \in M_{q^{\ast }}$ such that%
\begin{equation}
\liminf_{n\rightarrow \infty }|G_{n,\theta +\gamma /\sqrt{n},\sigma
}(t)-G_{\infty ,\theta +\gamma /\sqrt{n},\sigma }(t)|>0
\label{insert_condition}
\end{equation}%
holds. Note that (\ref{insert_condition}) is equivalent to%
\begin{equation}
\liminf_{n\rightarrow \infty }|G_{\infty ,\theta ,\sigma ,\gamma
}(t)-G_{\infty ,\theta +\gamma /\sqrt{n},\sigma }(t)|>0  \label{condition}
\end{equation}%
in light of Proposition \ref{p3.1}. To establish (\ref{condition}), we
proceed as follows: For each $\gamma \in M_{q^{\ast }}$ with $\gamma
_{q^{\ast }}\neq 0$, $G_{\infty ,\theta +\gamma /\sqrt{n},\sigma }(t)$ in (%
\ref{3.7}) reduces to $\Phi _{\infty ,q^{\ast }}(t)$ as is easily seen from (%
\ref{ab.1}) since $p_{0}(\theta +\gamma /\sqrt{n})=q^{\ast }$ which in turn
follows from $p_{0}(\theta )<q^{\ast }$ and $\gamma _{q^{\ast }}\neq 0$.
Furthermore, Lemma~\ref{lb.1} entails that $G_{\infty ,\theta ,\sigma
,\gamma }(t)$ is non-constant in $\gamma \in M_{q^{\ast }}\backslash
M_{q^{\ast }-1}$. But this shows that (\ref{condition}) must hold.

To prove part (\ref{p.ae1.b}), we write 
\begin{eqnarray*}
\left\vert \left\vert G_{n,\theta ,\sigma }-\Phi _{\infty ,P}\right\vert
\right\vert _{TV} &=&\left\vert \left\vert \sum_{p=\O }^{P}G_{n,\theta
,\sigma }(\cdot |p)\pi _{n,\theta ,\sigma }(p)-\Phi _{\infty ,P}(\cdot
)\right\vert \right\vert _{TV} \\
&\leq &\sum_{p=\O }^{P}\left\vert \left\vert G_{n,\theta ,\sigma }(\cdot
|p)-\Phi _{\infty ,P}(\cdot )\right\vert \right\vert _{TV}\pi _{n,\theta
,\sigma }(p),
\end{eqnarray*}%
where the conditional cdfs $G_{n,\theta ,\sigma }(\cdot |p)$ and the model
selection probabilities $\pi _{n,\theta ,\sigma }(p)$ have been introduced
after (\ref{3.5}). By the `uncorrelatedness' assumption, we have that $\Phi
_{\infty ,p}=\Phi _{\infty ,P}$ for all $p$ in the range $\O \leq p\leq P$;
cf. Remark \ref{rb.3}. We hence obtain%
\begin{equation}
\sup_{\theta \in \mathbf{R}^{P}}\sup_{\overset{\sigma \in \mathbf{R}}{\sigma
_{\ast }\leq \sigma \leq \sigma ^{\ast }}}\left\vert \left\vert G_{n,\theta
,\sigma }-\Phi _{\infty ,P}\right\vert \right\vert _{TV}\leq \sum_{p=\O %
}^{P}\sup_{\theta \in \mathbf{R}^{P}}\sup_{\overset{\sigma \in \mathbf{R}}{%
\sigma _{\ast }\leq \sigma \leq \sigma ^{\ast }}}\left\vert \left\vert
G_{n,\theta ,\sigma }(\cdot |p)-\Phi _{\infty ,p}(\cdot )\right\vert
\right\vert _{TV}\pi _{n,\theta ,\sigma }(p).  \label{sum}
\end{equation}%
Now for every $p$ with $\O \leq p\leq P$ and for every $\rho $, $0<\rho
<\infty $, we can write%
\begin{align}
& \sup_{\theta \in \mathbf{R}^{P}}\sup_{\overset{\sigma \in \mathbf{R}}{%
\sigma _{\ast }\leq \sigma \leq \sigma ^{\ast }}}\left\vert \left\vert
G_{n,\theta ,\sigma }(\cdot |p)-\Phi _{\infty ,p}(\cdot )\right\vert
\right\vert _{TV}\pi _{n,\theta ,\sigma }(p)  \notag \\
& \leq \max \left\{ \sup_{\overset{\theta \in \mathbf{R}^{P}}{\left\Vert
\theta \lbrack \lnot p]\right\Vert <\rho /\sqrt{n}}}\sup_{\overset{\sigma
\in \mathbf{R}}{\sigma _{\ast }\leq \sigma \leq \sigma ^{\ast }}}\left\vert
\left\vert G_{n,\theta ,\sigma }(\cdot |p)-\Phi _{\infty ,p}(\cdot
)\right\vert \right\vert _{TV}\,\ ,\sup_{\overset{\theta \in \mathbf{R}^{P}}{%
\left\Vert \theta \lbrack \lnot p]\right\Vert \geq \rho /\sqrt{n}}}\sup_{%
\overset{\sigma \in \mathbf{R}}{\sigma _{\ast }\leq \sigma \leq \sigma
^{\ast }}}\pi _{n,\theta ,\sigma }(p)\right\} .  \label{max}
\end{align}%
In case $p=P$, we use here the convention that the second term in the
maximum is absent and that the first supremum in the first term in the
maximum extends over all of $\mathbb{R}^{P}$. Letting first $n$ and then $%
\rho $ go to infinity in (\ref{max}), we may apply Lemmas C.2 and C.3 in
Leeb and P\"{o}tscher (2005b) to conclude that the l.h.s. of (\ref{max}),
and hence the l.h.s. of (\ref{sum}), goes to zero as $n\rightarrow \infty $.~%
{\hfill $\Box $}

\section{Proofs for Sections~\protect\ref{s3} to~\protect\ref{s5}\label{ac}}

In the proofs below it will be convenient to show the dependence of $\Phi
_{n,p}$ and $\Phi _{\infty ,p}$\ on $\sigma $ in the notation. Thus, in the
following we shall write $\Phi _{n,p,\sigma }$ and $\Phi _{\infty ,p,\sigma
} $, respectively, for the cdf of a $k$-variate Gaussian random vector with
mean zero and variance-covariance matrix $\sigma ^{2}A[p](X[p]^{\prime
}X[p]/n)^{-1}A[p]^{\prime }$ and $\sigma ^{2}A[p]Q[p:p]^{-1}A[p]^{\prime }$,
respectively. For convenience, let $\Phi _{n,0,\sigma }$ and $\Phi _{\infty
,0,\sigma }$ denote the cdf of point-mass at zero in $\mathbf{R}^{k}$.

The following lemma is elementary to prove, if we recall that $b_{n,p}z$
converges to $b_{\infty ,p}z$\ as $n\rightarrow \infty $ for every $z\in 
\func{Im}A[p]$, the column space of $A[p]$.

\begin{lemma}
\label{Lb.1}Suppose $p>\O $. Define $R_{n,p}(z,\sigma )=1-\Delta _{\sigma
\zeta _{n,p}}(b_{n,p}z,c_{p}\sigma \xi _{n,p})$ and\ $R_{\infty ,p}(z,\sigma
)=1-\Delta _{\sigma \zeta _{\infty ,p}}(b_{\infty ,p}z,c_{p}\sigma \xi
_{\infty ,p})$ for $z\in \func{Im}A[p]$, $0<\sigma <\infty $. Let $\sigma
^{(n)}$ converge to $\sigma $, $0<\sigma <\infty $. If $\zeta _{\infty
,p}\neq 0$, then $R_{n,p}(z,\sigma ^{(n)})$ converges to $R_{\infty
,p}(z,\sigma )$ for every $z\in \func{Im}A[p]$; if $\zeta _{\infty ,p}=0$,
then convergence holds for every $z\in \func{Im}A[p]$, except possibly for $%
z\in \func{Im}A[p]$ satisfying $\left\vert b_{\infty ,p}z\right\vert
=c_{p}\sigma \xi _{\infty ,p}$. [This exceptional subset of $\func{Im}A[p]$
has $rank(A[p])$-dimensional Lebesgue measure zero since $c_{p}\sigma \xi
_{\infty ,p}>0$.]
\end{lemma}

The following observation is useful in the proof of Proposition~\ref{p4.1.1}
below: Since the proposition depends on $Y$ only through its distribution
(cf. Remark~\ref{r6.0}), we may assume without loss of generality that the
errors in (\ref{2.1}) are given by $u_{t}=\sigma \varepsilon _{t}$, $t\in 
\mathbf{N}$, with i.i.d. $\varepsilon _{t}$ that are standard normal. In
particular, all random variables involved are then defined on the same
probability space.

\textbf{Proof of Proposition \ref{p4.1.1}:} Since $P_{n,\theta ,\sigma }(%
\bar{p}=p_{0}(\theta ))\rightarrow 1$ by consistency, we may replace $\max \{%
\bar{p},\O \}$ by $p_{\ast }=\max \{p_{0}(\theta ),\O \}$ in the formula for 
$\check{G}_{n}$\ for the remainder of the proof. Furthermore, since $\hat{%
\sigma}\rightarrow \sigma $ in $P_{n,\theta ,\sigma }$-probability, each
subsequence contains a further subsequence along which $\hat{\sigma}%
\rightarrow \sigma $\ almost surely (with respect to the probability measure
on the common probability space supporting all random variables involved),
and we restrict ourselves to such a further subsequence for the moment. In
particular, we write $\left\{ \hat{\sigma}\rightarrow \sigma \right\} $\ for
the event that $\hat{\sigma}$ converges to $\sigma $ along the subsequence
under consideration; clearly, the event $\left\{ \hat{\sigma}\rightarrow
\sigma \right\} $ has probability one. Also note that we can assume without
loss of generality that $\hat{\sigma}>0$ holds on this event (at least from
some data-dependent $n$ onwards), since $\sigma >0$ holds. But then
obviously $\prod_{q=p_{\ast }+1}^{P}\Delta _{\hat{\sigma}\xi _{n,q}}(0,c_{q}%
\hat{\sigma}\xi _{n,q})$ converges to $\prod_{q=p_{\ast }+1}^{P}\Delta
_{\sigma \xi _{\infty ,q}}(0,c_{q}\sigma \xi _{\infty ,q})$, and $\hat{\Phi}%
_{n,p_{\ast }}(t)$ converges to $\Phi _{\infty ,p_{\ast },\sigma }(t)$ in
total variation by Lemma A.3 of Leeb (2005) in case $p_{\ast }>0$, and
trivially so in case $p_{\ast }=0$. This proves that the first term in the
formula for $\check{G}_{n}$ converges to the corresponding term in the
formula for $G_{\infty ,\theta ,\sigma }$ in total variation.

Next, consider the term in $\check{G}_{n}$\ that carries the index $%
p>p_{\ast }$. By Lemma A.3 in Leeb (2005), $\hat{\Phi}_{n,p}=\Phi _{n,p,\hat{%
\sigma}}$ has a density $d\Phi _{n,p,\hat{\sigma}}/d\Phi _{\infty ,p,\sigma
} $ with respect to $\Phi _{\infty ,p,\sigma }$, which converges to $1$
except on a set that has measure zero under $\Phi _{\infty ,p,\sigma }$. By
Scheff\'{e}'s Lemma (Billingsley (1995), Theorem 16.12), $d\Phi _{n,p,\hat{%
\sigma}}/d\Phi _{\infty ,p,\sigma }$ converges to $1$ also in the $%
L^{1}(\Phi _{\infty ,p,\sigma })$-sense. By Lemma \ref{Lb.1}, $R_{n,p}(z,%
\hat{\sigma})$ converges to $R_{\infty ,p}(z,\sigma )$ except possibly on a
set that has measure zero under $\Phi _{\infty ,p,\sigma }$. (Recall that $%
\Phi _{\infty ,p,\sigma }$ is concentrated on $\func{Im}A[p]$ and is not
degenerate there.) Observing that $\left\vert R_{n,p}(z,\hat{\sigma}%
)\right\vert $ is uniformly bounded by $1$, we obtain that $R_{n,p}(z,\hat{%
\sigma})$ converges to $R_{\infty ,p}(z,\sigma )$ also in the $L^{1}(\Phi
_{\infty ,p,\sigma })$-sense. Hence, 
\begin{align}
& \left\Vert R_{n,p}(z,\hat{\sigma})\frac{d\Phi _{n,p,\hat{\sigma}}}{d\Phi
_{\infty ,p,\sigma }}(z)-R_{\infty ,p}(z,\sigma )\right\Vert  \notag \\
& \leq \left\Vert R_{n,p}(z,\hat{\sigma})\frac{d\Phi _{n,p,\hat{\sigma}}}{%
d\Phi _{\infty ,p,\sigma }}(z)-R_{n,p}(z,\hat{\sigma})\right\Vert
+\left\Vert R_{n,p}(z,\hat{\sigma})-R_{\infty ,p}(z,\sigma )\right\Vert
\label{ac.0} \\
& \leq \left\Vert \frac{d\Phi _{n,p,\hat{\sigma}}}{d\Phi _{\infty ,p,\sigma }%
}(z)-1\right\Vert +\left\Vert R_{n,p}(z,\hat{\sigma})-R_{\infty ,p}(z,\sigma
)\right\Vert \overset{n\rightarrow \infty }{\longrightarrow }0  \notag
\end{align}%
where $\left\Vert \cdot \right\Vert $ denotes the $L^{1}(\Phi _{\infty
,p,\sigma })$-norm. Since $\prod_{q=p+1}^{P}\Delta _{\hat{\sigma}\xi
_{n,q}}(0,c_{q}\hat{\sigma}\xi _{n,q})$ obviously converges to $%
\prod_{q=p+1}^{P}\Delta _{\sigma \xi _{\infty ,q}}(0,c_{q}\sigma \xi
_{\infty ,q})$, the relation (\ref{ac.0}) shows that the term in $\check{G}%
_{n}$\ carrying the index $p$ converges to the corresponding term in $%
G_{\infty ,\theta ,\sigma }$ in the total variation sense. This proves (\ref%
{4.3}) along the subsequence under consideration. However, since any
subsequence contains such a further subsequence, this establishes (\ref{4.3}%
). Since $G_{n,\theta ,\sigma }$ converges to $G_{\infty ,\theta ,\sigma }$
in total variation by Proposition \ref{p3.1}, the claim in (\ref{4.2}) also
follows.{\hfill $\Box $}

Before we prove the main result we observe that the total variation distance
between $P_{n,\theta ,\sigma }$ and $P_{n,\vartheta ,\sigma }$ satisfies $%
\left\vert \left\vert P_{n,\theta ,\sigma }-P_{n,\vartheta ,\sigma
}\right\vert \right\vert _{TV}\leq 2\Phi (\left\Vert \theta -\vartheta
\right\Vert \lambda _{\max }^{1/2}(X^{\prime }X)/2\sigma )-1$; furthermore,
if $\theta ^{(n)}$ and $\vartheta ^{(n)}$ satisfy $\left\Vert \theta
^{(n)}-\vartheta ^{(n)}\right\Vert =O(n^{-1/2})$, the sequence $%
P_{n,\vartheta ^{(n)},\sigma }$ is contiguous with respect to the sequence $%
P_{n,\theta ^{(n)},\sigma }$ (and vice versa). This follows exactly in the
same way as Lemma~A.1 in Leeb and P\"{o}tscher (2006a).

\textbf{Proof of Theorem \ref{t5.1}: }We first prove (\ref{5.2}) and (\ref%
{5.3}). For this purpose we make use of Lemma~3.1 in Leeb and P\"{o}tscher
(2006a) with $\alpha =\theta \in M_{q^{\ast }-1}$, $B=M_{q^{\ast }}$, $%
B_{n}=\{\vartheta \in M_{q^{\ast }}:\left\Vert \vartheta -\theta \right\Vert
<\rho _{0}n^{-1/2}\}$, $\beta =\vartheta ,$ $\varphi _{n}(\beta
)=G_{n,\vartheta ,\sigma }(t)$, $\widehat{\varphi }_{n}=\hat{G}_{n}(t)$,
where $\rho _{0}$, $0<\rho _{0}<\infty $, will be chosen shortly (and $%
\sigma $ is held fixed). The contiguity assumption of this lemma (as well as
the mutual contiguity assumption used in the corrigendum to Leeb and P\"{o}%
tscher (2006a)) is satisfied in view of the preparatory remark above. It
hence remains only to show that there exists a value of $\rho _{0}$, $0<\rho
_{0}<\infty $, such that $\delta ^{\ast }$ in Lemma~3.1 of Leeb and P\"{o}%
tscher (2006a) (which represents the limit inferior of the oscillation of $%
\varphi _{n}(\cdot )$ over $B_{n}$) is positive. Applying Lemma~3.5(i) of
Leeb and P\"{o}tscher (2006a) with $\zeta _{n}=\rho _{0}n^{-1/2}$ and the
set $G_{0}$ equal to the set $G$, it remains, in light of Proposition~\ref%
{p3.1}, to show that there exists a $\rho _{0}$, $0<\rho _{0}<\infty $, such
that $G_{\infty ,\theta ,\sigma ,\gamma }(t)$ as a function of $\gamma $ is
non-constant on the set $\{\gamma \in M_{q^{\ast }}:\left\Vert \gamma
\right\Vert <\rho _{0}\}$. In view of Lemma~3.1 of Leeb and P\"{o}tscher
(2006a), the corresponding $\delta _{0}$ can then be chosen as any positive
number less than one-half of the oscillation of $G_{\infty ,\theta ,\sigma
,\gamma }(t)$ over this set. That such a $\rho _{0}$ indeed exists follows
now from Lemma~\ref{lb.1} in Appendix \ref{aa}, where it is also shown that $%
\rho _{0}$ and $\delta _{0}$ can be chosen such that they depend only on $%
t,Q,A,\sigma ,$ and $c_{p}$ for $\O <{p}\leq P$. This completes the proof of
(\ref{5.2}) and (\ref{5.3}).

To prove (\ref{5.4}) we use Corollary~3.4 in Leeb and P\"{o}tscher (2006a)
with the same identification of notation as above, with $\zeta _{n}=\rho
_{0}n^{-1/2}$, and with $V=M_{q^{\ast }}$ (viewed as a vector space
isomorphic to $\mathbf{R}^{q^{\ast }}$). The asymptotic uniform
equicontinuity condition in that corollary is then satisfied in view of $%
\left\vert \left\vert P_{n,\theta ,\sigma }-P_{n,\vartheta ,\sigma
}\right\vert \right\vert _{TV}\leq 2\Phi (\left\Vert \theta -\vartheta
\right\Vert \lambda _{\max }^{1/2}(X^{\prime }X)/2\sigma )-1$. Given that
the positivity of $\delta ^{\ast }$ has already be established in the
previous paragraph, applying Corollary~3.4(i) in Leeb and P\"{o}tscher
(2006a) then establishes (\ref{5.4}). {\hfill $\Box $}

\textbf{Proof of Remark \ref{r5.1}: }The proof is similar to the proof of (%
\ref{5.4}) just given, except for using Corollary~3.4(ii) and Lemma 3.5(ii)
in Leeb and P\"{o}tscher (2006a) instead of Corollary 3.4(i) and Lemma
3.5(i) from that paper. Furthermore, Lemma \ref{lb.2} in Appendix \ref{aa}
instead of Lemma \ref{lb.1} is used.{\hfill $\Box $}

\textbf{Proof of Proposition \ref{p5.2}: }In view of Proposition~\ref{p.ae1}(%
\ref{p.ae1.b}) and the fact that $\hat{\Phi}_{n,P}(\cdot )=\Phi _{n,P,\hat{%
\sigma}}(\cdot )$ holds (in case $\hat{\sigma}>0$), it suffices to show that%
\begin{equation}
\sup_{\overset{\sigma \in \mathbf{R}}{\sigma _{\ast }\leq \sigma \leq \sigma
^{\ast }}}\,\left\vert \left\vert \Phi _{n,P,\sigma }(\cdot )-\Phi _{\infty
,P,\sigma }(\cdot )\right\vert \right\vert _{TV}\overset{n\rightarrow \infty 
}{\longrightarrow }0  \label{C.4'}
\end{equation}%
\begin{equation}
\sup_{\overset{\sigma \in \mathbf{R}}{\sigma _{\ast }\leq \sigma \leq \sigma
^{\ast }}}\,P_{n,\theta ,\sigma }\left( \left\vert \left\vert \Phi _{n,P,%
\hat{\sigma}}(\cdot )-\Phi _{n,P,\sigma }(\cdot )\right\vert \right\vert
_{TV}\,>\,\delta \right) \overset{n\rightarrow \infty }{\longrightarrow }0
\label{C.4}
\end{equation}%
hold for each $\delta >0$, and for any constants $\sigma _{\ast }$ and $%
\sigma ^{\ast }$ satisfying $0<\sigma _{\ast }\leq \sigma ^{\ast }<\infty $.
[Note that the probability in (\ref{C.4}) does in fact not depend on $\theta 
$.] But this has already been established in the proof of Proposition 4.3 of
Leeb and P\"{o}tscher (2005b).{\hfill $\Box $}

\section{Proofs for Section \protect\ref{s99}\label{ad}}

\textbf{Proof of Theorem \ref{t99.2}: }After rearranging the elements of $%
\theta $ (and hence the regressors) if necessary and then correspondingly
rearranging the rows of the matrix $A$, we may assume without loss of
generality that $\mathfrak{r}_{\ast }=(1,\ldots ,1,0)$, and hence that $i(%
\mathfrak{r}_{\ast })=P$. That is, $M_{\mathfrak{r}_{\ast }}=M_{P-1}$ and $%
M_{\mathfrak{r}_{full}}=M_{P}$. Furthermore, note that after this
arrangement $C_{\infty }^{(P)}\neq 0$. Let $\hat{p}$\ be the model selection
procedure introduced in Section \ref{s2} with $\O =P-1$, $c_{P}=c$, and $c_{%
\O }=0$. Let $\tilde{\theta}$ be the corresponding post-model-selection
estimator and let $G_{n,\theta ,\sigma }(t)$ be as defined in Section \ref%
{s3}. Condition (\ref{99.1}) now implies: For every $\theta \in M_{P-1}$
which has exactly $P-1$ non-zero coordinates 
\begin{equation}
\lim_{n\rightarrow \infty }P_{n,\theta ,\sigma }\left( \{\hat{\mathfrak{r}}=%
\mathfrak{r}_{full}\}\blacktriangle \{\hat{p}=P\}\right) =\lim_{n\rightarrow
\infty }P_{n,\theta ,\sigma }\left( \{\hat{\mathfrak{r}}=\mathfrak{r}_{\ast
}\}\blacktriangle \{\hat{p}=P-1\}\right) =0  \label{E.1}
\end{equation}%
holds for every $0<\sigma <\infty $. Since the sequences $P_{n,\vartheta
^{(n)},\sigma }\ $and $P_{n,\theta ,\sigma }$ are contiguous for $\vartheta
^{(n)}$ satisfying $\left\Vert \theta -\vartheta ^{(n)}\right\Vert
=O(n^{-1/2})$ as remarked prior to the proof of Theorem \ref{t5.1} in
Appendix \ref{ac}, it follows that condition (\ref{E.1}) continues to hold
with $P_{n,\vartheta ^{(n)},\sigma }$ replacing $P_{n,\theta ,\sigma }$.
This implies that for every sequence of positive real numbers $s_{n}$ with $%
s_{n}=O(n^{-1/2})$, for every $\sigma $, $0<\sigma <\infty $, and for every $%
\theta \in M_{P-1}$ which has exactly $P-1$ non-zero coordinates 
\begin{equation}
\sup_{\overset{\vartheta \in \mathbf{R}^{P}}{||\vartheta -\theta ||<s_{n}}%
}\left\Vert K_{n,\vartheta ,\sigma }-G_{n,\vartheta ,\sigma }\right\Vert
_{TV}\rightarrow 0  \label{E.4}
\end{equation}%
holds as $n\rightarrow \infty $. From (\ref{E.4}) we conclude that the limit
of $K_{n,\theta +\gamma /\sqrt{n},\sigma }$\ (with respect to total
variation distance) exists and coincides with $G_{\infty ,\theta ,\sigma
,\gamma }$. Repeating the proof of Theorem~\ref{t5.1} with $q^{\ast }=P$,
with $K_{n,\vartheta ,\sigma }(t)$ replacing $G_{n,\vartheta ,\sigma }(t)$,
and with $\hat{K}_{n}(t)$ replacing $\hat{G}_{n}(t)$ gives the desired
result.{\hfill $\Box $}

\section{References}

\quad\ Ahmed, S.~E. \& A.~K. Basu (2000): Least squares, preliminary test
and Stein-type estimation in general vector AR(p) models. \emph{Statistica
Neerlandica \ }54, 47--66.

Bauer, P., P\"{o}tscher, B.~M. \& P. Hackl (1988): Model selection by
multiple test procedures. \emph{Statistics \ }19, 39--44.

Billingsley, P. (1995): \emph{Probability and Measure,} (3rd ed.).\emph{\ }%
Wiley.

Brownstone, D. (1990): Bootstrapping improved estimators for linear
regression models. \emph{Journal of Econometrics \ }44, 171--187.

Danilov, D.~L. \& J.~R. Magnus (2004): On the harm that ignoring pre-testing
can cause. \emph{Journal of Econometrics \ }122, 27--46.

Dijkstra, T.~K. \& J.~H. Veldkamp (1988): `Data-driven selection of
regressors and the bootstrap'. \emph{Lecture Notes in Economics and
Mathematical Systems \ }307, 17--38.

Duki\'{c}, V.~M. \& E.~A Pe\~{n}a (2002): Variance estimation in a model
with gaussian submodel. \emph{Journal of the American Statistical Association%
} \ 100, 296-309.

Freedman, D.~A., Navidi, W. \& S.~C. Peters (1988): `On the impact of
variable selection in fitting regression equations'. \emph{Lecture Notes in
Economics and Mathematical Systems \ }307, 1--16.

Hansen, P.~R. (2003): Regression analysis with many specifications: a
bootstrap method for robust inference. Working Paper, Department of
Economics, Brown University.

Hjort, N.~L. \& G. Claeskens (2003): Frequentist model average estimators. 
\emph{Journal of the American Statistical Association} \ 98, 879--899.

Kabaila, P. (1995): The effect of model selection on confidence regions and
prediction regions. \emph{Econometric Theory} {\ 11}, 537--549.

Kapetanios, G. (2001): Incorporating lag order selection uncertainty in
parameter inference for AR models. \emph{Economics Letters } 72, 137--144.

Kilian, L. (1998): Accounting for lag order uncertainty in autoregressions:
the endogenous lag order bootstrap algorithm. \emph{Journal of Time Series
Analysis \ }19, 531--548.

Knight, K. (1999): Epi-convergence in distribution and stochastic
equi-semicontinuity. Working Paper, Department of Statistics, University of
Toronto.

Kulperger, R.~J. \& S.~E.\ Ahmed (1992): A bootstrap theorem for a
preliminary test estimator. \emph{Communications in Statistics: Theory and
Methods }\ 21, 2071--2082.

Leeb, H. (2002): On a differential equation with advanced and retarded
arguments. \emph{Communications on Applied Nonlinear Analysis} 9, 77--86.

Leeb, H. (2005): The distribution of a linear predictor after model
selection: conditional finite-sample distributions and asymptotic
approximations. \emph{Journal of Statistical Planning and Inference \ }134,
64--89.

Leeb, H. (2006): The distribution of a linear predictor after model
selection: unconditional finite-sample distributions and asymptotic
approximations. \emph{IMS Lecture Notes-Monograph Series \ }49, 291--311.

Leeb, H. \& B.~M. P\"{o}tscher (2003): The finite-sample distribution of
post-model-selection estimators and uniform versus nonuniform
approximations. \emph{Econometric Theory} {\ 19}, 100--142.

Leeb, H. \& B.~M. P\"{o}tscher (2005a): Model selection and inference: facts
and fiction. \emph{Econometric Theory} {\ 21}, 21--59.

Leeb, H. \& B.~M. P\"{o}tscher (2005b): Can one estimate the conditional
distribution of post-model-selection estimators? Working Paper, Department
of Statistics, University of Vienna.

Leeb, H. \& B.~M. P\"{o}tscher (2006a): Performance limits for estimators of
the risk or distribution of shrinkage-type estimators, and some general
lower risk bound results. \emph{Econometric Theory \ }22, 69-97.
(Corrigendum. \emph{Econometric Theory}, forthcoming.)

Leeb, H. \& B.~M. P\"{o}tscher (2006b): Can one estimate the conditional
distribution of post-model-selection estimators? \emph{Annals of Statistics}
\ 34, 2554-2591.

Lehmann, E.~L. \& G. Casella (1998): \emph{Theory of Point Estimation}, 2nd
Edition, Springer Texts in Statistics. Springer-Verlag.

Nickl, R. (2003): \emph{Asymptotic Distribution Theory of
Post-Model-Selection Maximum Likelihood Estimators. }Masters Thesis,
Department of Statistics, University of Vienna.

P\"{o}tscher, B.~M. (1991): Effects of model selection on inference. \emph{%
Econometric Theory} {\ 7}, 163--185.

P\"{o}tscher, B.~M\textsc{.} (1995): {Comment} on {`The} effect of model
selection on confidence regions and prediction regions' by P. Kabaila. \emph{%
Econometric Theory} {\ 11}, 550--559.

P\"{o}tscher, B.~M. \& A.~J. Novak (1998): The distribution of estimators
after model selection: large and small sample results. \emph{Journal of
Statistical Computation and Simulation} {\ 60}, 19--56.

Rao, C.~R. \& Y. Wu (2001): `On model selection,' \emph{IMS Lecture
Notes-Monograph Series \ }38, 1--57.

Robinson, G.~K. (1979): Conditional properties of statistical procedures, 
\emph{Annals of Statistics } 7, 742--755.

Sen, P.~K. (1979): Asymptotic properties of maximum likelihood estimators
based on conditional specification. \emph{Annals of Statistics \ }7,
1019--1033.

Sen P.~K. \& A.~K.~M.~E. Saleh (1987): On preliminary test and shrinkage
M-estimation in linear models. \emph{Annals of Statistics \ }15, 1580--1592.

van der Vaart, A.~W. (1998): \emph{Asymptotic Statistics. }Cambridge
University Press.

\bibliographystyle{econometrica}
\bibliography{lit}

\end{document}